\documentclass[12pt]{article}
\usepackage[bottom]{footmisc}

\usepackage[utf8]{inputenc}
\usepackage[english]{babel}
\usepackage[T1]{fontenc}
\usepackage{hyperref}
\hypersetup{
   colorlinks,
    linkcolor={red!50!black},
    citecolor={green!50!black},
    urlcolor={blue!80!black}
}
\usepackage{array}
\usepackage{multicol}
\usepackage{bbm}
\usepackage{amscd}
\usepackage{amsmath}
\usepackage{nccmath}
\usepackage{mathtools, bm}
\usepackage{amssymb, bm}
\usepackage{amsthm}
\usepackage{stmaryrd}
\usepackage[unq]{unique}
\usepackage{thmtools}
\usepackage{amsfonts}
\usepackage{dsfont}
\usepackage{graphics}
\usepackage{tikz}
\usepackage{tikz-cd}
\usepackage{geometry}
\usepackage{enumerate}

\usepackage{setspace}
\usepackage{parskip}
\usepackage{soul} 
\usepackage{url}

\usepackage[all]{xy}

\declaretheoremstyle[
  headfont=\itshape, qed = \qedsymbol
]{prv}
\declaretheoremstyle[
    headfont=\bfseries, bodyfont=\itshape
]{theoreme}
\declaretheoremstyle[
  headfont=\bfseries,
]{definition}

\declaretheorem[style=theoreme, name=Theorem]{thm}
\declaretheorem[style=theoreme, name=Proposition,sibling=thm]{prop}
\declaretheorem[style=theoreme, name=Lemma,sibling=thm]{lemme}
\declaretheorem[style=theoreme, name=Corollary,sibling=thm]{coro}
\declaretheorem[style=definition, name=Remark,sibling=thm]{rmq}

\declaretheorem[style=definition, name=Definition,sibling=thm]{definition}
\declaretheorem[style=prv, name=Proof, numbered = no]{prv}

\DeclareMathOperator{\Gal}{Gal}

\newcommand{\mbb}{\mathbb}
\newcommand{\mbf}{\mathbf}
\newcommand{\mc}{\mathcal}
\newcommand{\op}{\operatorname}
\newcommand{\abs}[1]{\left| #1 \right|}

\newcommand{\N}{\mathbb{N}}
\newcommand{\Z}{\mathbb{Z}}
\newcommand{\Q}{\mathbb{Q}}
\newcommand{\R}{\mathbb{R}}
\newcommand{\C}{\mathbb{C}}
\newcommand{\K}{\mathbb{K}}

\DeclareMathOperator{\was}{\mathbf{Was}}

\DeclareMathOperator{\Sin}{\mbf{Sin}}

\DeclareMathOperator{\frob}{\op{Frob}}

\allowdisplaybreaks

\begin{document}

\baselineskip=17pt

\title{Bases for some modules of cyclotomic units}

\author{Rafik SOUANEF}

\date{}

\maketitle


\renewcommand{\thefootnote}{}

\footnote{2020 \emph{Mathematics Subject Classification}: Primary 11R27, 11R18; Secondary 11R23.}

\footnote{\emph{Key words and phrases}: Cyclotomic units, real cyclotomic fields, totally deployed fields, bases, basis, generators, circular units}

\renewcommand{\thefootnote}{\arabic{footnote}}
\setcounter{footnote}{0}

\vspace{-2cm}

\begin{abstract}
    Let $\was(\mbb{K})$ denote the group of Washington's cyclotomic units of any abelian number field $\mbb{K}$. If $\mbb{K}$ coincides with its genus field in the narrow sense, we give a $\Lambda$-basis of $\lim\limits_{\xleftarrow{}} \was(\K_k^{+})$ where $(\K_k)_{k \geqslant 0}$ denotes the cyclotomic $\Z_p$-tower of $\K$ and $\Lambda$ denotes the Iwasawa's algebra. This results from a $\Z [1/2]$-basis of $\was(\mbb{K}) \otimes_{\Z} \Z [1/2]$ that we give under the same hypothesis.
\end{abstract}

\section{Introduction}

Let $\K$ be an abelian number field of conductor $n$. Let $\zeta_n = \op{exp}(2i \pi / n)$. Let $\mbf{Z}(\K)$ denote the group formed by the roots of unity of $\K$ and let $\mbf{E}(\K)$ denote the group of units (of the ring of integers) of $\K$. Let $\mc{C}_n$ be the Galois module generated by $\mbf{Z}(\Q(\zeta_n))$ and the family $(1-\zeta_d)_{d > 1, d \mid n}$. Let $\was(\mbb{K}) = \mbf{E}(\mbb{K}) \cap \mc{C}_n $ denote the group of Washington's cyclotomic units of $\K$ \cite[before Lemma 8.1]{was}. Let $\Sin(\mbb{K})$ denote the group of Sinnott's cyclotomic units of $\K$, that is the intersection of $\mbf{E}(\mbb{K})$ with the Galois module generated by $\mbf{Z}(\K)$ and the family $(\op{N}_{\mbb{Q}(\zeta_d) / \mbb{K} \cap \mbb{Q}(\zeta_d)} ( 1- \zeta_d))_{d > 1, d \mid n}$ \cite[before Lemma 4.1]{Sinnott1980}. The abelian group $\Sin(\K)$ comes with explicit generators \cite[before section 2]{Kuc97} and a formula for $[\mbf{E}(\K) : \Sin(\K)]$ in which the class number $h(\K^{+})$ appears \cite[Theorem 4.1]{Sinnott1980}. On the other hand, no such results are known for $\was(\K)$ in general. Explicit $\Z$-bases of $\was(\K)$ have been obtained for some fields that coincide with their genus field \cite[Corollary 4.3]{KuceraBase}, \cite[Proposition 2, Remark 4]{werl}. Also, the difference between these two types of cyclotomic units is not clear in general. However, a criterion was given in \cite[Proposition 3.6, Theorem 2.2]{bel} for these two types of units to coincide at infinity. At infinity, the major difference between these units is that Washington's units always form a free module over the Iwasawa's algebra \cite[before Proposition 3.6]{bel}. The reader may see \cite{KuceraRevue} for an exposition on the different definitions of cyclotomic units. In particular, finding bases of $\was(\K)$ and $\Sin(\K)$ at both finite and infinite level is still an open question and, in this paper, we provide some answers to these questions with Theorems \ref{thmBinf1}, \ref{thmBinf2} and \ref{item:11}. We also give an upper bound on the cardinality of the not so well understood quotient $\was(\K)/\Sin(\K)$ in Corollary \ref{diffwassin}. We will focus on totally deployed (abelian number) fields, that is to say fields of the form $\mbb{K} = \K (1) \cdots \K (r) $ with $\K (i) \subset \mbb{Q}(\zeta_{p_i^{e_i}})$ for some distinct prime numbers $p_i$ and some integers $e_i$. In other words, we are interested in abelian number fields that coincide with their genus field in the narrow sense.

Theorems \ref{thmBinf1} and \ref{thmBinf2} consist in giving a $\Lambda$-basis of $\lim\limits_{\xleftarrow{}} \was(\K_k^{+})$, where $(\K_k)_{k \geqslant 0}$ denotes the cyclotomic $\Z_p$-tower of any totally deployed field $\K$, $\Lambda=\lim\limits_{\xleftarrow{}} \Z_p[\Gal(\K_k/\K)]$ denotes the Iwasawa's algebra and $p$ denotes any odd prime integer. The need of having two theorems is explained by a distinction of cases that we make for technical reasons.

Finding this basis at infinity mostly relies on Theorem \ref{item:11} in which we give a $\Z [1/2]$-basis of $\was(\K) \otimes \Z[1/2]$ assuming $\K$ is totally deployed. Let $B(\Q(\zeta_n))$ denote the $\Z$-basis of $\was(\Q(\zeta_n))$ given by \cite[Corollary 4.3]{KuceraBase} (for convenience, we recall this result in Theorem \ref{theoremkucbase}). The main idea in the proof of Theorem \ref{item:11} is to show that we have a $\Z [1/2]$-basis by proving that the elements we consider generate a direct factor of $\was(\mbb{Q}(\zeta_n)) \otimes \Z [1/2]$. In order to use this idea, we have to make quite technical computations that result from an algorithm that allows us to compute the decomposition of any element of $\was(\Q(\zeta_n))$ in the basis $B(\Q(\zeta_n))$; this algorithm is partially described with Lemmas \ref{keyf1} and \ref{keyf2}. This idea has also been used in \cite[Theorem 2.1]{kimryu} and \cite[Proposition 2]{werl} to get $\Z$-bases of $\was(\K)$ in a setup that requires direct computations.

Divisibility relations on class numbers arise from the basis given by Theorem \ref{item:11} (see Corollary \ref{item:19}) and results of Sinnott on $[\mbf{E}(\K) : \Sin(\K)]$.

\section{Preliminaries}

Let $\N$ denote the set of all natural integers and let $\N^{*}=\mbb{N} \setminus \{ 0 \}$ be the set of all positive integers.

\subsection{On units} \label{onunits}

Let $\K$ be an abelian number field. Let $A$ be a $\Z[\Gal(\K/\Q)]$-module. Observe that the complex conjugation induces an element of $\Gal(\mbb{K}/\mbb{Q})$. Let $A^{+}$ denote the Galois submodule of $A$ that consists of all the elements of $A$ on which the complex conjugation acts trivially. Later on, we will focus on multiplicative structures. For any $x \in A$ and for any $u \in \mbb{Z}[\Gal(\mbb{K}/\mbb{Q})]$, we denote the image of $x$ under $u$ by $ u(x)$ or $x^u$.

Let $\zeta_n = \op{exp}(2i \pi / n)$ for any $n \in \mbb{N}^{*}$. From now on, let $n \geqslant 2$ satisfy $n \not\equiv 2 \mod 4$. If $p$ is a prime number, let $v_p(k)$ denote the $p$-valuation of any integer $k$. 

For the rest of the article, let $\mbb{K}$ be an abelian number field of conductor $n$ (this explains the condition $n \not\equiv 2 \mod 4$). Let $ n = \prod_{j=1}^r p_j ^{e_j}$ with $p_j$ being a prime number and $e_j \in \mbb{N}$ for any $j \in \llbracket 1,r \rrbracket$. Let $q_j = p_j^{e_j}$ for any $j \in \llbracket 1,r \rrbracket$.
We say $\mbb{K}$ is totally deployed when $\Gal(\mbb{K}/\mbb{Q})$ is the direct product of its inertia subgroups (see the introduction of \cite{gras}). As we supposed $\K/\Q$ is abelian, the field $\K$ is totally deployed if and only if, for any $j \in \llbracket 1,r \rrbracket$, there is a number field $\K (j) \subset \mbb{Q}(\zeta_{q_j})$  such that
\[ \mbb{K} = \K (1) \cdots \K (r) .\]
Let $\mbf{E}(\mbb{K})$ be the group of units (of the ring of integers $\mc{O}_{\mbb{K}}$) of $\mbb{K}$. Recall that if $n$ is not a prime power, then $1-\zeta_n \in \mbf{E}(\Q(\zeta_n))$ (see \cite[Proposition 2.8]{was}). If $n$ is a prime power, then $1-\zeta_n$ is no longer a unit but $(1-\zeta_n^\sigma)/(1-\zeta_n)$ is a unit for any $\sigma \in \Gal(\mbb{Q}(\zeta_n) / \mbb{Q} )$ (see \cite[Lemma 1.3]{was}). Let $\mbf{Z}(\K)$ denote the group of roots of unity of $\K$.

\begin{definition}
Let $\mc{C}_n$ be the Galois module generated by $\mbf{Z}(\mbb{Q}(\zeta_n))$ and by the $1-\zeta_d$'s with $d \mid n, d > 1$. Let $\was(\mbb{K}) = \mbf{E}(\mbb{K}) \cap \mc{C}_n $. Let $\Sin(\mbb{K}) $ be the intersection of $\mbf{E}(\mbb{K})$ with the Galois module generated by $\mbf{Z}(\K)$ and the family $(\op{N}_{\mbb{Q}(\zeta_d) / \mbb{K} \cap \mbb{Q}(\zeta_d)} ( 1- \zeta_d))_{d > 1} .$
\end{definition}



When the situation makes it clear, we will omit writing $\mbb{K}$. For example, we will write $\was$ instead of writing $\was(\mbb{K})$.

It is known that cyclotomic units satisfy the following relations (see \cite{solomon}, Lemma 2.1):
\begin{equation} \label{conj}
    1-\zeta_n^a = -\zeta_n^a (1- \zeta_n^{-a})
\end{equation}
\begin{equation} \label{normrelcyc}\op{N}_{\mbb{Q}(\zeta_n) / \mbb{Q}(\zeta_d)} ( 1- \zeta_n) =  \left(\prod\limits_{\substack{p | n \\ p \nmid d}} (1-\frob(p)^{-1}) \right) (1-\zeta_d)\end{equation}
where $d \mid n$ is such that $d > 1$, the integers $p$ are prime and $\frob(p)$ denotes the Frobenius of $\mbb{Q}(\zeta_d)$ that is defined by $\zeta_d \mapsto \zeta_d^p$. We will refer to this second relation as 'norm relation'. We will call this relation 'norm relation along $\sigma_i$' (we will define $\sigma_i$ later) to mean that we consider this norm relation with $d=n/q_i$. These relations are our main tools to prove Lemmas \ref{keyf1} and \ref{keyf2}, in which we partially explain how to decompose elements of $\was(\Q(\zeta_n))$ in the basis of $\was(\Q(\zeta_n))$ given by \cite[Corollary 4.3]{KuceraBase} (for convenience, we recall this corollary of \cite{KuceraBase} in Theorem \ref{theoremkucbase}).

We recall a property of Hasse's unit index.

\begin{prop}[\cite{was}, Theorem 4.12, Corollary 4.13.] \label{item:10}
    We have
    \[ [\mathbf{E}: \mathbf{Z}\mathbf{E}^{+}] \in \{ 1;2 \} .\]
Moreover, if $\mbb{K}=\mbb{Q}(\zeta_n)$, this index is $1$ if and only if $n$ is a prime power.
\end{prop}

Remind Dirichlet's units theorem: the abelian group $\mbf{E}(\mbb{K})$ is finitely generated, its torsion part is $\mbf{Z}(\mbb{K})$ and it has rank $r_1+r_2-1$ (where $r_1$ is the number of real embeddings of $\K$ and $r_2$ is half of the number of complex embeddings of $\K$).
It is known that both $\was(\mbb{K})$ and $\Sin(\mbb{K})$ have finite index in $\mbf{E}(\mbb{K})$, that is they both have the same rank as $\mbf{E}(\mbb{K})$ (see \cite[Theorem 4.1]{Sinnott1980}). This allows us to use a simple strategy to prove that some family $F$ is a basis of $\was(\K)$: if $F$ has cardinality $r_1+r_2-1$ and generates $\was(\K)$, then $F$ is a basis of $\was(\K)$.

We now introduce some of the notation that we will use to work with bases of $\was$ (most of this notation comes from \cite{gold}, \cite{werl} and \cite{KuceraBase}).

For any $j \in \llbracket 1,r \rrbracket$, let $J_j$ denote the complex conjugation considered as an element of $\Gal(\mbb{Q}(\zeta_{q_j})/\mbb{Q})$. If $p_j$ is odd, let $\sigma_j$ be a generator of $\Gal(\mbb{Q}(\zeta_{q_j})/\mbb{Q})$. If $p_j=2$, let $\sigma_j$ be such that $\Gal(\mbb{Q}(\zeta_{q_j})/\mbb{Q})$ is generated by $\sigma_j$ and $J_j$ (so that $\Gal(\mbb{Q}(\zeta_{q_j})/\mbb{Q})$ is the direct product of $\langle J_j \rangle$ and $ \langle  \sigma_j \rangle$).

From now on, for any $j \in \llbracket 1,r \rrbracket$, see the elements of $\Gal(\mbb{Q}(\zeta_{q_j})/\mbb{Q})$ as elements of $\Gal(\mbb{Q}(\zeta_{n})/\mbb{Q})$ by letting them act trivially on $\mbb{Q}(\zeta_{n/q_j})$. Let $J= J_1 \cdots J_r$ be the complex conjugation considered as an element of $\Gal(\mbb{Q}(\zeta_{n})/\mbb{Q})$.

For any $u \in \Gal(\Q(\zeta_n)/\Q)$, define up to a sign
\[\xi_{q_j,u} = \sqrt{\zeta_{q_j}^{1-u}}\frac{1-\zeta_{q_j}^{u}}{1-\zeta_{q_j}} \in \was(\mbb{Q}(\zeta_{q_j})^{+}). \]

We now introduce the notation we will use to define the basis of $\was(\Q(\zeta_n))$ given by \cite[Corollary 4.3]{KuceraBase} (see Theorem \ref{theoremkucbase}). The reader must be aware that, when one talks about bases, it is common to talk about $\was(\mbb{K})$ instead of the free abelian group $\was(\mbb{K}) / \mbf{Z}(\mbb{K})$.


\begin{definition}[\cite{KuceraBase}, Lemma 1.1] \label{defRi}
    For any $i \in \llbracket 1,r \rrbracket$, a set $\mc{R}_i$ is defined  in the following way. If $p_i \neq 2$, let $z \in \Gal(\Q(\zeta_{q_i}) / \Q)$ be such that $z$ generates the $2$-Sylow of $\Gal(\Q(\zeta_{q_i}) / \Q)$ and let $H$ be the non $2$-part of $\Gal(\Q(\zeta_{q_i}) / \Q)$ (that is $H$ is the product of the Sylow $l$-subgroups of $\Gal(\Q(\zeta_{q_i}) / \Q)$ for $l$ running over the set of odd prime numbers). Let $a \in \N$ be such that $z^{2^{a}}=J_i$. Let
    \[
    \mc{R}_i = \{ z^{k} h : 0 \leqslant k < 2^a, h \in H\}.
    \]
 If $p_i=2$, let $\mc{R}_i = \langle \sigma_i \rangle$.
\end{definition}

\begin{rmq}
    For any $i \in \llbracket 1,r \rrbracket$, the set $\mc{R}_i$ is a set of representatives of $\Gal(\mbb{Q}(\zeta_{q_i})/\Q)$ modulo $\langle J_i \rangle$ and we have $1 \in \mc{R}_i$.
\end{rmq}

\begin{definition}
    Let $\Omega \subset \llbracket 1,r \rrbracket$ be a non-empty set.
    If $\Omega=\{ i \} \subset \llbracket 1,r \rrbracket$ for some integer $i$, let $X_{\Omega}$ denote $\mc{R}_i \setminus \{ 1 \}$.
    Otherwise, let $s \geqslant 2$ and $i_1 < \dots < i_s$ be such that $\Omega=\{ i_1, \dots, i_s \}$. Let $X_{\Omega}$ be the set of products of the form $u_1 \cdots u_k$ with $k \in \llbracket 1,s \rrbracket$, satisfying $u_k \in \mc{R}_{i_k} \setminus \{ 1 \}$ and 
    \[
    \forall j \in \llbracket 1, k-1 \rrbracket, \quad u_j \in \Gal(\mbb{Q}(\zeta_{q_{i_j}})/\Q) \setminus \{ J_{i_j} \}. 
    \]
    If $\abs{\Omega}$ is even, then add $1$ to $X_{\Omega}$.
\end{definition}

\begin{definition} \label{defcn}
    For any non-empty set $\Omega \subset \llbracket 1,r \rrbracket $, let $n_{\Omega}=\prod_{i \in \Omega} q_i$, let $\zeta_{\Omega}=\zeta_{n_{\Omega}}$ and let
    $
    b_{n_{\Omega}} = b_{\Omega} = 1-\zeta_{\Omega}
    $.
    If $\Omega=\{ i \}$ for some $i \in \llbracket 1,r \rrbracket$, let
    \[
    B_{\Omega} = \left\{ \xi_{q_i,u} : u \in X_{\Omega} \right\}.
    \]
    If $\abs{\Omega} \geqslant 2$, let
    \[
    B_{\Omega} = \left\{ b_{\Omega}^{u} : u \in X_{\Omega} \right\}.
    \]
\end{definition}

\begin{thm}[\cite{KuceraBase}, Corollary 4.3] \label{theoremkucbase}
    The family $B = \cup_{\Omega} B_{\Omega}$ where $\Omega$ runs over the set of all non-empty subsets of $\llbracket 1,r \rrbracket$ is a basis of $\was(\Q(\zeta_n))$.
\end{thm}

In the rest of the text, we may write $B(\Q(\zeta_d))$ (or $B_{\Omega}(\Q(\zeta_d))$) to talk about the basis that is given by Theorem \ref{theoremkucbase} for $\Q(\zeta_d)$, with $d$ being any positive integer that satisfies $d \not\equiv 2 \mod 4$. We now make some remarks on this theorem.

\begin{definition} \label{deflevel}
Let $\Omega$ be a non-empty subset of $\llbracket 1,r \rrbracket$. We say $\Omega$ or $n_{\Omega}$ is the level of any element of $B_{\Omega}$.
\end{definition}

We have $B(\Q(\zeta_{d})) \subset B(\Q(\zeta_{d'})) $ for any $d \mid d'$ such that $(d'/d) \wedge d = 1$ so that any element of $\was(\Q(\zeta_d))$ decomposes in $B(\Q(\zeta_{d'}))$ with terms whose level is lower than or equal to $d$.

Theorem \ref{theoremkucbase} comes with an algorithm to compute the decomposition of any element of $\was(\mbb{Q}(\zeta_n))$ in the basis $B$ (see the proof of Lemma 3.2 in \cite{KuceraBase}). In order to prove Theorem \ref{item:11}, we partially explain this algorithm and make additional observations on those decompositions in Lemmas \ref{keyf1} and \ref{keyf2}. These two lemmas are essential to prove Theorem \ref{item:11}. If $n$ is not a prime power, we will explain what amounts to decomposing in $\cup_{\abs{\Omega} \geqslant 2} B_{\Omega}(\Q(\zeta_n))$ any
\[
b_n^{u_1 \cdots u_r} \in \was(\Q(\zeta_n)) / \langle B_{\{i\}}(\Q(\zeta_n)), i \in \llbracket 1,r \rrbracket \rangle
\]
given $u_1 \cdots u_r \in \Gal(\Q(\zeta_n)/\Q) \setminus \langle J_1, \dots, J_r \rangle$ such that $u_i \in \Gal(\Q(\zeta_{q_i})/\Q)$ for any $i \in \llbracket 1,r \rrbracket$. This algorithm is partially presented here because we will neither be concerned with the decomposition of any $b_n^u$ with $u \in  \langle J_1, \dots, J_r \rangle$ nor will we need to make a whole lemma on the elements of the $\was(\Q(\zeta_{q_i}))$'s (for $i$ running over $ \llbracket 1,r \rrbracket$). This algorithm works by induction on $r$ (that is the number of distinct prime factors of $n$). However, in the following, we will not mention any induction hypothesis on $r$ as we will be interested only in the elements of $B_{\llbracket 1,r \rrbracket}$ that appear in the decomposition of those $b_n^{u_1 \cdots u_r}$. In other words, any element of $\was(\Q(\zeta_d))$ that appear in the following - with $d$ being a strict divisor of $n$ - may be decomposed by induction and this is enough information for us.

Let
\begin{gather*}
    V(u_1, \dots, u_r) = \{ i \in \llbracket 1,r \rrbracket : u_i = J_i \} \\
    W(u_1, \dots, u_r) = \{ i \in \llbracket 1,r \rrbracket : u_i \neq 1 \}
\end{gather*}
for any $u_i \in \Gal(\Q(\zeta_{q_i})/ \Q)$, $i \in \llbracket 1,r \rrbracket$. For Lemmas \ref{keyf1} and \ref{keyf2}, suppose $r \geqslant 2$. Lemma \ref{keyf1} states that we may suppose $u_i \neq J_i$ for any $i \in \llbracket 1,r \rrbracket$, which then allows us to be in the setup that is needed for Lemma \ref{keyf2}.

\begin{lemme} \label{keyf1}
For any $i \in \llbracket 1,r \rrbracket$, let $u_i \in \Gal(\Q(\zeta_{q_i})/\Q)$ and suppose that we have $u_1 \cdots u_r \not \in \langle J_1, \dots, J_r \rangle$. Let $b= b_n^{u_1 \cdots u_r}$.
The unit $b$ is a product of elements whose level is lower than $n$ and terms of the form $b_n^{\pm w_1 \cdots w_r}$ with $w_i \in \Gal(\Q(\zeta_{q_i})/\Q) \setminus \{ J_i \}$ for any $i \in \llbracket 1,r \rrbracket$ and $w_i = u_i$ for any $i \not\in V(u_1, \dots, u_r)$.
\end{lemme}

\begin{prv}
    We will work by induction on $\abs{V(u_1, \dots, u_r)}$. Indeed, if $\abs{V(u_1, \dots, u_r)}=0$, there is nothing to do. Suppose $\abs{V(u_1, \dots, u_r)}>=1$ and assume Lemma \ref{keyf1} holds for any $w_1, \dots, w_r$ such that $V(w_1,\dots,w_r) < V(u_1,\dots,u_r)$. There is $i \in V(u_1, \dots, u_r)$ and the norm relation along $\sigma_i$ gives
\[
    b_n^{u_1 \cdots u_{r}} = b_{n/q_i}^{(1-\frob(p_i)^{-1}) u_1 \cdots u_{r}} \prod\limits_{w_i \in \Gal(\mbb{Q}(\zeta_{q_i})/\Q) \setminus \{ J_i \}} b_n^{- u_1 \cdots u_{i-1} w_i u_{i+1} \cdots u_{r}}.
\]
Observe $b_{n/q_i}^{(1-\frob(p_i)^{-1})u_1 \cdots u_{r}} \in \was(\Q(\zeta_{n/q_i}))$ decomposes in $B(\Q(\zeta_{n/q_i}))$ (so that it decomposes with elements of $B(\Q(\zeta_n))$ whose level is lower than $n$) and the elements of the form $b_n^{- u_1 \cdots u_{i-1} w_i u_{i+1} \cdots u_{r}}$ that appear above are such that 
\[
V(u_1, \dots u_{i-1}, w_i, u_{i+1}, \dots, u_{r}) \subsetneq V(u_1, \dots, u_r).
\]
Then, using the induction hypothesis on those $b_n^{u_1 \cdots u_{i-1} w_i u_{i+1} \cdots u_{r}}$ concludes.
\end{prv}

\begin{lemme} \label{keyf2}
    For any $i \in \llbracket 1,r \rrbracket$, let $u_i \in \Gal(\Q(\zeta_{q_i})/ \Q) \setminus \{ J_i \}$ and assume that we have $u_1 \cdots u_r \neq 1$. Let $b= b_n^{u_1 \cdots u_r}$ and let $u=u_1 \cdots u_r$. Let $\max W(u_1, \dots, u_r)$ be the length of $u$ and $b$; denote it by $M(u)$.
     We have the following two points
    \begin{enumerate}[i)]
    \item \label{prempoint} if some $b_n^{w_1 \cdots w_k} \in B$ appears in the decomposition of $b$, then we have $k \geqslant M(u)$.
    \item \label{deuzpoint} Moreover, assume $u_{M(u)} \in \mc{R}_{M(u)}$ or $u_i \neq 1$ for any $i < M(u)$; if some $b_n^{w_1 \cdots w_{M(u)}} \in B$ appears in the decomposition of $b$, then we have
    \[
    w_1 \cdots w_{M(u)} = \left\{
        \begin{array}{ll}
            u & \mbox{ if } u_{M(u)} \in \mc{R}_{M(u)} \\
            J_1 u_1 \cdots J_{M(u)}u_{M(u)} & \mbox{ otherwise.}
        \end{array}
    \right.
    \]
    Also, $b_n^{w_1 \cdots w_{M(u)}}$ appears with exponent $1$ in the first case and exponent $(-1)^{r-M(u)}$ in the second.
    \end{enumerate}
\end{lemme}

\begin{prv}
    We will show  that by backward induction on $M(u)$.
  
\underline{The base case} corresponds to the case $M(u)=r$. Then, assume $M(u)=r$. If $u_r \in \mc{R}_r$, then we have $u \in X_{\llbracket 1,r \rrbracket}$ so that we have $b \in B$ and there is nothing to do. Now, assume $u_r \not\in \mc{R}_r$ (so that we actually have $u_r \in J_r \mc{R}_r \setminus \{ J_r \}$). Equation \eqref{conj} gives, modulo roots of unity
\[
b=b_n^{J_1 u_1 \cdots J_r u_r}
\]
and we have $J_r u_r \in \mc{R}_r \setminus \{1\}$. It suffices to call Lemma \ref{keyf1} (for $J_1 u_1 \cdots J_r u_r$) to get \ref{prempoint}). Moreover, suppose $u_i \neq 1$ for any $i < r$, then we have $J_1 u_1 \cdots J_r u_r \in X_{\llbracket 1, r \rrbracket}$ and this concludes the proof of the base case.

\underline{For the induction step}, suppose there is $k \in \llbracket 1,r-1 \rrbracket$ such that Lemma \ref{keyf2} holds whenever $M(u) \geqslant k+1$. Assume $M(u) = k$. If $u_{k} \in \mc{R}_k$, then we have $u \in X_{\llbracket 1,r \rrbracket}$ so there is nothing to do. Now, assume $u_k \not\in \mc{R}_k$, that is $u_k \in J_k \mc{R}_k \setminus \{ J_k \}$. We will use norm relations consecutively to conclude. More precisely, the norm relation along $\sigma_{k+1}$ gives, modulo elements whose level is lower than $n$
\[
b=b_n^{-u_1 \cdots u_k J_{k+1}} \prod\limits_{u_{k+1} \in \Gal(\Q(\zeta_{q_{k+1}})/\Q) \setminus \{1,J_{k+1} \}} b_n^{-u_1 \cdots u_{k+1}}.
\]
Those $u_1 \cdots u_{k+1}$ are such that $M(u_1 \cdots u_{k+1}) \geqslant k+1$ and we have $u_i \neq J_i$ for any $i \in \llbracket 1, k+1 \rrbracket$. Then, the induction hypothesis shows these elements decompose in $B$ with lower level elements and elements of $B_{\llbracket 1,r \rrbracket}(\Q(\zeta_n))$ that have length greater than $k$. We can repeat this process by using the norm relation along $\sigma_{k+2}$ on $b_n^{-u_1 \cdots u_k J_{k+1}}$ etc. so that we will prove by induction on $l \in \llbracket k+1, r \rrbracket$ that we have
\begin{equation} \label{eqdecompo}
b=b_n^{(-1)^{l-k} u_1 \cdots u_k J_{k+1} \cdots J_l}
\end{equation}
modulo $\mbf{Z}(\Q(\zeta_n))$, modulo elements whose level is lower than $n$ and modulo elements of $B_{\llbracket 1,r \rrbracket}(\Q(\zeta_n))$ that have length greater than $k$. Indeed, we just proved the base case $l=k+1$. Assume Equation \eqref{eqdecompo} holds for some $l$ and let us prove it holds for $l+1$. The norm relation along $\sigma_{l+1}$ gives
\begin{align*}
b_n^{(-1)^{l-k} u_1 \cdots u_k J_{k+1} \cdots J_l} =& b_n^{(-1)^{l+1-k}u_1 \cdots u_k J_{k+1} \cdots J_{l+1}} \\ 
& \prod\limits_{\substack{u_{l+1} \in \Gal(\Q(\zeta_{q_{l+1}})/ \Q) \\ u_{l+1} \neq 1, J_{l+1}}} b_n^{(-1)^{l+1-k}u_1 \cdots u_k J_{k+1} \cdots J_{l} u_{l+1}}
\end{align*}
modulo elements whose level is lower than $n$.
It follows from Lemma \ref{keyf1} that any such $b_n^{(-1)^{l+1-k}u_1 \cdots u_k J_{k+1} \cdots J_{l} u_{l+1}}$ is a product of elements whose level is lower than $n$ and some $b_n^{\pm u_1 \cdots u_k w_{k+1} \cdots w_{l} u_{l+1}}$ with $w_{k+1} \neq J_{k+1}, \dots, w_l \neq J_l$. Then, the induction hypothesis on $M$ shows that any of those $b_n^{\pm u_1 \cdots u_k w_{k+1} \cdots w_{l} u_{l+1}}$ decomposes in $B$ with elements whose level is lower than $n$ and elements of $B_{\llbracket 1,r \rrbracket}(\Q(\zeta_n))$ that have length greater than $l$. This concludes the proof by induction on $l$.
Hence, taking $l=r$ in Equation \eqref{eqdecompo} before using Equation \eqref{conj} gives
\begin{equation} \label{fineqdecompo}
b=b_n^{(-1)^{r-k} J_1 u_1 \cdots J_k u_k}
\end{equation}
modulo $\mbf{Z}(\Q(\zeta_n))$, modulo elements whose level is lower than $n$ and modulo elements of $B_{\llbracket 1,r \rrbracket}(\Q(\zeta_n))$ that have length greater than $k$.
Then, Lemma \ref{keyf1} gives \ref{prempoint}). If $u_i \neq 1$ for any $i < k$, we have $J_1 u_1 \cdots J_k u_k \in X_{\llbracket 1,r \rrbracket}$ and Equation \eqref{fineqdecompo} concludes the proof by induction on $M$.
\end{prv}

\subsection{On the convolution product} \label{sectionconvprod}

Through this section, we recall - in the needed context only - some facts that are stated in a more general context in \cite{rota}, \cite{greene} and that deal with Möbius functions.

Let $E$ be a finite set and let $\mc{P}(E)$ denote the powerset of $E$. Define $\mc{F}(E)$ as the set of functions
\[ \begin{array}{lccc} f: & \mc{P}(E) & \longrightarrow & \mbb{C} \end{array}. \]
This set has a law of addition and a convolution product defined in the following way
\[\forall f,g \in \mc{F}(E), \ \forall \Omega \subset E, \quad  f * g (\Omega) = \sum\limits_{X \subset \Omega} f(X) g(\Omega \setminus X). \]

One can show $(\mc{F}(E),+,*)$ is a ring whose identity element is the function that maps $\emptyset$ to $1$ and any subset $\Omega \neq \emptyset$ to $0$.

Denote by $\mbf{1}$ the element of $\mc{F}(E)$ that maps any $\Omega \subset E$ to $1$.
One can show $\mbf{1}$ is a unit and we let $\mu$ denote its inverse. We have (see \cite[Equation 3.3]{greene})
\[ \forall \Omega \subset E, \quad \mu(\Omega) = (-1)^{|\Omega|} .\]
In particular, we have the following theorem.

\begin{thm}[\cite{rota}, Proposition 2] \label{thmconvprod}
    Let $f,g \in \mc{F}(E)$. We have
    \[  \forall \Omega \subset E, \; \sum\limits_{X \subset \Omega} f(X) = g(\Omega) \Longleftrightarrow{}  \forall \Omega \subset E, \; f(\Omega) = \sum\limits_{X \subset \Omega} (-1)^{|\Omega| - |X|} g(X). \]
\end{thm}

Later, we will use this convolution product with $E = \llbracket 1 , r \rrbracket$ to prove that the family $C(\K)$ - that is considered in Theorem \ref{item:11} - has cardinality $\op{rank}(\was(\K))$ (see the beginning of the proof of Theorem \ref{item:11}).

\section{Totally deployed fields} 

Let $\was_2(\K) = \was(\K) \otimes_{\Z} \Z [1/2]$. In this section, we give a $\Z [1/2]$-basis of $\was_2(\mbb{K})$ (see Theorem \ref{item:11}) assuming $\mbb{K}$ is a totally deployed abelian number field. In particular, we will have a family that is a $\Z_p$-basis of $\was(\K) \otimes \Z_p$ for any prime integer $p > 2$ and this will be our starting point in the construction of $\Lambda$-bases in section \ref{sectioninfini}. For now, we suppose $\mbb{K}$ is a totally deployed abelian number field of conductor $n$ and we write
\[ \mbb{K}= \K (1) \cdots \K (r)\]
with $\K (i) \subset \mbb{Q}(\zeta_{q_i})$ for any $i \in \llbracket 1,r \rrbracket$. To simplify the proof of Theorem \ref{item:11}, if there is $i$ such that $p_i=2$ and $\K (i)$ is imaginary, suppose $i=r$.

To construct our basis, we will consider a family $C(\K)$ that generates a direct factor of $\was_2(\mbb{Q}(\zeta_n))$ and that is made of $ r_1+r_2-1$ elements of $\K$. It is not hard to see that this property makes $C(\K)$ generate $\was_2(\mbb{K})$: this is what Lemma \ref{item:18} shows. Then, the family $C(\K)$ is a basis of $\was_2(\mbb{K})$. This idea has already been used in \cite[Proposition 2]{werl}, \cite[Theorem 2.1]{kimryu}. Actually, in order to prove \cite[Proposition 2]{werl}, the author proves our Lemma \ref{item:18} for abelian groups.

\begin{lemme} \label{item:18}
    Let $R$ be a principal ideal domain. Let $M_1 \subset M_2 \subset M_3$ be free $R$-modules. Assume $M_1$ and $M_2$ have same rank. Suppose $M_1$ is a direct factor of $M_3$. Then, we have $M_1=M_2$.
\end{lemme}

\begin{proof}
    We simply adapt the proof of \cite[Proposition 2]{werl}. There is $r \in R \setminus \{ 0 \}$ such that $rM_2 \subset M_1$ since $M_1$ and $M_2$ have same rank and $R$ is a principal ideal domain. Let $m_2 \in M_2$. Let $M_4$ be such that $M_3 = M_1 \oplus M_4$. There are $m_1 \in M_1$, $m_4 \in M_4$ such that $m_2=m_1+m_4$. We have
    \[
    \underbrace{r m_2}_{\in M_1} = \underbrace{r m_1}_{\in M_1} + \underbrace{r m_4}_{\in M_4}
    \]
    hence the definition of $M_4$ implies $m_4=0$, that is to say $m_2=m_1 \in M_1$.
\end{proof}

\subsection{Notation} \label{sectionnotCK}

Recall $\mc{R}_i$ is the set of representatives of $\Gal(\Q(\zeta_{q_i})/ \Q)$ modulo $J_i$ given by \cite[Lemma 1.1]{KuceraBase} (see Definition \ref{defRi}). We now introduce the notation we will use to state Theorem \ref{item:11}.

Up to reordering, there is $t$ such that $\K (1),\dots,\K (t-1)$ are real and $\K (t),\dots,\K (r)$ are imaginary.

For any $i \in \llbracket 1, t-1 \rrbracket$, let $(\mc{R}_{i,1}(\K), \mc{T}_{i}(\K))$ be such that $\mc{R}_{i,1}(\K)$ is a set of representatives of $\Gal(\Q(\zeta_{q_i}) / \Q) / \Gal(\Q(\zeta_{q_i})/ \K (i))$ with $1 \in \mc{R}_{i,1}(\K)$ and $\mc{T}_{i}(\K)$ is a set of representatives of $\Gal(\Q(\zeta_{q_i})/ \K (i)) / \langle J_i \rangle $ such that $\mc{T}_i \cdot \mc{R}_{i,1}(\K) \subset \mc{R}_i$. \newline
For instance, we can construct $\mc{R}_{i,1}(\K)$ and $ \mc{T}_{i}(\K)$ as follows. First, if $p_i=2$ then that construction is clear as $\langle J_i \rangle$ is a direct factor of $\Gal(\Q(\zeta_{q_i}) / \Q)$. Then, suppose $p_i$ is odd. Recall $z$ denotes a generator of the $2$-Sylow of $\Gal(\Q(\zeta_{q_i}) / \Q)$ (see Definition \ref{defRi}) and let $m \in \N$ be minimal with respect to $z^{2^{m}} \in \Gal(\Q(\zeta_{q_i}) / \K (i))$. Let $a \in \N$ be such that $z^{2^{a}} = J_i$. Let 
    \begin{align*}
         &\mc{T}_i(\K) = \{ z^{k2^m} h : k \in \llbracket 0 , 2^{a-m} \llbracket , h \in \Gal(\Q(\zeta_{q_i}) / \K (i)) \mbox{ has odd order} \} \\
         &\mc{R}_{i,1}(\K) = \{ z^k h : 0 \leqslant k < 2^m \mbox{ and } h \in H(\K) \}
    \end{align*}
where $H(\K)$ denotes any set of representatives of the non $2$-part of $\Gal(\Q(\zeta_{q_i}) / \Q)$ modulo $\Gal(\Q(\zeta_{q_i}) / \K (i))$ that lies in the non $2$-part of $\Gal(\Q(\zeta_{q_i}) / \Q)$. Now, swap $1$ with $J_i$ in $\mc{R}_{i,1}(\K)$.

For any $i \in \llbracket t, r \rrbracket$, let $\mc{R}_{i,1}(\K)$ be a set of representatives of $\Gal(\Q(\zeta_{q_i}) / \Q)$ modulo $\Gal(\Q(\zeta_{q_i}) / \K (i))$ with $1, J_i \in \mc{R}_{i,1}(\K)$. If $p_i \neq 2$, observe $\Gal(\Q(\zeta_{q_i})/ \K (i))$ acts on $\mc{R}_i$ by multiplication. Then, let $\mc{R}_{i,2}(\K)$ be a set of representatives of $\mc{R}_i$ modulo $\Gal(\Q(\zeta_{q_i})/ \K (i))$ with $1 \in \mc{R}_{i,2}(\K)$. In particular $\mc{R}_{i,2}(\K)$ is a set of representatives of $\Gal(\Q(\zeta_{q_i}) / \Q)$ modulo $\langle J_i, \Gal(\Q(\zeta_{q_i}) / \K (i)) \rangle$.
If $p_i=2$, observe $\Gal(\Q(\zeta_{q_i})/ \K (i))$ still acts on $\mc{R}_i$ and define $\mc{R}_{i,2}(\K)$ as before (that action is given by a transport of structure through the canonical bijection $\mc{R}_i \simeq \Gal(\Q(\zeta_{q_i})/\Q) / \langle J_i \rangle$ as $\Gal(\Q(\zeta_{q_i})/ \K (i))$ acts on this last quotient by multiplication). The set $\mc{R}_{i,2}(\K)$ is still a set of representatives of $\Gal(\Q(\zeta_{q_i}) / \Q)$ modulo $\langle J_i, \Gal(\Q(\zeta_{q_i}) / \K (i)) \rangle$ but we can no longer assume $\Gal(\Q(\zeta_{q_i})/ \K (i)) \cdot \mc{R}_i \subset \mc{R}_i$.

Let
 $\mbb{L}_n = \mbb{Q}(\zeta_{q_1})^{+} \cdots \mbb{Q}(\zeta_{q_r})^{+} .$
If $r \geqslant 2$, there is a root of unity $\eta \in \mbb{Q}(\zeta_n)$ \cite[2-ii)]{werl} such that
$\eta_n= \eta \op{N}_{\mbb{Q}(\zeta_n) / \mbb{Q}(\zeta_{q_1})^{+} \cdots \mbb{Q}(\zeta_{q_{r-1}})^{+} \mbb{Q}(\zeta_{q_r}) } (1-\zeta_n) \in \mbb{Q}(\zeta_{q_1})^{+} \cdots \mbb{Q}(\zeta_{q_r})^{+}$
and $\eta_n^2 = \op{N}_{\mbb{Q}(\zeta_n) / \mbb{Q}(\zeta_{q_1})^{+} \cdots \mbb{Q}(\zeta_{q_{r-1}})^{+} \mbb{Q}(\zeta_{q_r})^{+} } (1-\zeta_n)$.
For any field $\mbb{L} \subset \mbb{L}_n$ whose conductor is $n$, let
$ e_{\mbb{L}} = \op{N}_{\mbb{L}_n / \mbb{L}} (\eta_n) \in \was(\mbb{L}).$

For any non-empty ${\Omega} = \{ i_1,\dots,i_s \} \subset \llbracket 1 , r \rrbracket $, let
\begin{gather*}
    \mbb{K}_{\Omega} = \K (i_1) \cdots \K (i_s), \quad \Omega_{\mbb{R}} = \Omega \cap \llbracket 1,t-1 \rrbracket, \quad \Omega_{\mbb{C}} = \Omega \cap \llbracket t,r \rrbracket \\
    c_{n_{\Omega}}(\K) = c_{\Omega}(\K) = \renewcommand{\arraystretch}{1.5}\left\{
    \begin{array}{*3{>{\displaystyle}l}}
       N_{\Q(\zeta_{\Omega})^{+} / \K_{\Omega}^{+}}(\xi_{q_i,\sigma_i}) & \mbox{if } \Omega = \{i \} \mbox{ for some } i \in \llbracket 1,r \rrbracket \\
        N_{\Q(\zeta_{\Omega}) / \K_{\Omega}} (1-\zeta_{\Omega}) & \mbox{if } \abs{\Omega_{\C}} \geqslant 1 \\
        e_{\K_{\Omega}} & \mbox{otherwise.}
    \end{array}
\right. 
\end{gather*} 

For any ${\Omega} = \{ i \} \subset \llbracket 1 , r \rrbracket $, let $Y_{\Omega}(\K)$ be $\mc{R}_{i,1}(\K) \setminus \{ J_i \} $ if $i < t$ and let $Y_{\Omega}(\K)$ be $\mc{R}_{i,2}(\K) \setminus \{ 1 \} $ otherwise. Let $t_{\Omega}=2$.
  
For any ${\Omega} = \{ i_1,\dots,i_s \} \subset \llbracket 1 , r \rrbracket $ with $s \geqslant 2$, such that $i_1 < \dots < i_s$ and $\mbb{K} (i_{s})$ is real (that is $\mbb{K}_{\Omega}$ decomposes with real fields only), let
$Y_{\Omega}(\K)$ be the set of all $u_1 \cdots u_s$ such that $u_j \in \mc{R}_{i_j,1}(\K) \setminus \{ J_{i_j} \}$ for any $j \in \llbracket 1 ,s \rrbracket$. Let $t_{\Omega}=s+1$.

For any ${\Omega} = \{ i_1,\dots,i_s \} \subset \llbracket 1 , r \rrbracket $ with $s \geqslant 2$, such that $i_1 < \dots < i_s$ and $\K (i_s)$ is imaginary (that is $\mbb{K}_{\Omega}$ decomposes with at least one imaginary field), let $t_{\Omega}$ be the integer such that $\K (i_1),\dots,\mbb{K} (i_{t_{\Omega}-1})$ are real and $ \mbb{K} (i_{t_{\Omega}}), \dots, \K (i_s)$ are imaginary.
Let $Y_{\Omega}(\K)$ be the set of products of the form $u_1 \cdots u_k$ with $k \in \llbracket t_{\Omega},s \rrbracket$, satisfying $u_k \in \mc{R}_{i_k,2}(\K) \setminus \{ 1 \}$ and $u_j \in \mc{R}_{i_j,1}(\K) \setminus \{ J_{i_j} \}$ for any $j \in \llbracket 1,k-1 \rrbracket$.
 If $\abs{\Omega_{\C}}$ is even, then add to $Y_{\Omega}(\K)$ all the products of the form $u_1 \cdots u_{t_{\Omega}-1}$ with $u_j \in \mc{R}_{i_j,1}(\K) \setminus \{ J_{i_j} \}$ for any $j \in \llbracket 1 , t_{\Omega}-1 \rrbracket$ (if $t_{\Omega}=1$, understand that we add $1$ to $Y_{\Omega}(\K)$).

For any non-empty $\Omega \subset \llbracket 1,r \rrbracket$, let
    \[
    C_{\Omega}(\K) = \left\{ c_{\Omega}(\K)^u : u \in Y_{\Omega}(\K) \right\}.
    \]
Let $C(\mbb{K})=\cup_{\Omega} C_{\Omega}(\mbb{K})$ where ${\Omega}$ runs over the set of all non-empty subsets of $\llbracket 1,r \rrbracket $.
\begin{definition}
    If $\Omega \subset \llbracket 1,r \rrbracket$ is such that $\abs{\Omega_{\C}} \neq 0 $ is even and $\Omega_{\R} \neq \emptyset$, we will call 'problematic terms' those $u_1 \cdots u_{t_{\Omega}-1} \in Y_{\Omega}(\K)$ or their corresponding elements in $C_{\Omega}(\K)$.
\end{definition}
    
To understand the proof of Theorem \ref{item:11} more easily, we now highlight that we have
\[
\forall u,u' \in Y_{\llbracket 1,r \rrbracket}(\K), \quad u \neq u' \implies u \neq u' \mod \Gal(\Q(\zeta_n)/\K).
\]

To help the reader understand what could the family $C(\K)$ be, we give an example with the following diagram. Assume $r=3$, $p_1=23$, $p_2=11$, $p_3 = 19$, $n=p_1 p_2 p_3$ $\K=\K_1 \K_2 \K_3$ with $\K_1=\Q(\zeta_{13})^{+}$ being real, $\K_2=\Q(\zeta_{11})$ being imaginary, $\K_3 = \Q(\zeta_{19})$. In the following, $\sigma_i$ denotes a generator of $\Gal(\Q(\zeta_{p_i})/\Q)$. Observe we have $4 \nmid p_i-1$ for any $i$ so that $\sigma_i^{2}$ generates the non $2$-part of $\Gal(\Q(\zeta_{p_i})/\Q)$ and $J_i$ generates the Sylow $2$-subgroup of $\Gal(\Q(\zeta_{p_i})/\Q)$ - all that makes the construction of $C(\K)$ easier. In some sense, this is the easiest example that has not been considered in the literature before (the case with one real field and one imaginary field results of the case with two real fields) and covers different parity of $\abs{\Omega_{\C}}$ as $\Omega$ runs over the set of all non-empty subsets of $\llbracket 1,3 \rrbracket$. We have $[\K : \Q ] = 1980$ and $\op{rank}_{\Z}(\was(\K)) = 989$.

\begin{tikzpicture}
    \node (n) at (6.5,14) [draw,thick,minimum width=14cm,minimum height=3cm, label=$C_{\llbracket 1,3 \rrbracket}(\K)$] {} ;
    \path (1.75,14.6) node[align=center]{$N_{\Q(\zeta_n)/\K}(1-\zeta_n)^{\sigma_{1}^{2u_1}}$ \\ $u_1 \in \llbracket 1, 10 \rrbracket$} ;
    \path (6,14.3) node[align=center]{$N_{\Q(\zeta_n)/\K}(1-\zeta_n)^{\sigma_{1}^{2u_1}\sigma_{2}^{2u_2}}$ \\ $u_1 \in \llbracket 1, 10 \rrbracket$ \\ $u_2 \in \llbracket 1, 4 \rrbracket$} ;
    \path (10.75,14) node[align=center]{$N_{\Q(\zeta_n)/\K}(1-\zeta_n)^{\sigma_{1}^{2u_1}\sigma_{2}^{u_2}\sigma_{3}^{2u_3}}$ \\ $u_1 \in \llbracket 1, 10 \rrbracket$ \\ $u_2 \in \llbracket 0, 10 \rrbracket \setminus \{ 5 \} $  \\ $u_3 \in \llbracket 1, 8 \rrbracket$ } ;
    \node (q1q2) at (1,10) [draw,thick,minimum width=5.5cm,minimum height=2.5cm, label=$C_{\{1,2 \}}(\K)$] {} ;
    \path (1,10) node[align=center]{$N_{\Q(\zeta_{253})/\K_1\K_2}(1-\zeta_{253})^{\sigma_{1}^{2u_1} \sigma_{2}^{2u_2}}$ \\ $u_1 \in \llbracket 1, 10 \rrbracket$ \\ $u_2 \in \llbracket 1, 4 \rrbracket$ } ;
    \draw (n) -- (q1q2) ;
    \node (q1q3) at (6.5,10) [draw,thick,minimum width=5.5cm,minimum height=2.5cm,label={[xshift=2.5em]$C_{\{1,3 \}}(\K)$}] {} ;
    \path (6.5,10) node[align=center]{$N_{\Q(\zeta_{437})/\K_1\K_3}(1-\zeta_{437})^{\sigma_{1}^{2u_1} \sigma_{3}^{2u_2}}$ \\ $u_1 \in \llbracket 1, 10 \rrbracket$ \\ $u_2 \in \llbracket 1, 8 \rrbracket$ } ;
    \draw (n) -- (q1q3) ;
    \node (q2q3) at (12,9.2) [draw,thick,minimum width=5.5cm,minimum height=4.5cm,label=$C_{\{2,3 \}}(\K)$] {} ;
    \path (12,10) node[align=center]{$N_{\Q(\zeta_{209})/\K_2\K_3}(1-\zeta_{209})^{\sigma_{2}^{u_1} \sigma_{3}^{2u_2}}$ \\ $u_1 \in \llbracket 0, 10 \rrbracket \setminus \{ 5 \}$ \\ $u_2 \in \llbracket 1, 8 \rrbracket$ } ;
    \path (12,8) node[align=center]{$N_{\Q(\zeta_{209})/\K_2\K_3}(1-\zeta_{209})^{\sigma_{2}^{2u_1}}$ \\ $u_1 \in \llbracket 0, 4 \rrbracket$ } ;
    \draw (n) -- (q2q3);
    \node (q1) at (1,4) [draw,thick,minimum width=4cm,minimum height=2cm,label={[xshift=2em]above left:$C_{\{1 \}}(\K)$}] {} ;
    \path (1,4) node[align=center]{$N_{\Q(\zeta_{23})^{+}/\K_1}(\xi_{23,\sigma_1})^{\sigma_{1}^{2u_1}}$ \\ $u_1 \in \llbracket 1, 10 \rrbracket$};
    \node (q2) at (6.5,4) [draw,thick,minimum width=4cm,minimum height=2cm, label=$C_{\{2 \}}(\K)$] {} ;
    \path (6.5,4) node[align=center]{$N_{\Q(\zeta_{11})/\K_2}(\xi_{1,\sigma_2})^{\sigma_{2}^{2u_1}}$ \\ $u_1 \in \llbracket 1, 4 \rrbracket$} ;
    \node (q3) at (12,4) [draw,thick,minimum width=4cm,minimum height=2cm, label={[xshift=-2em]above right:$C_{\{3 \}}(\K)$}] {} ;
    \path (12,4) node[align=center]{$N_{\Q(\zeta_{19})/\K_3}(\xi_{19,\sigma_3})^{\sigma_{3}^{2u_1}}$ \\ $u_1 \in \llbracket 1, 8 \rrbracket$};
    \draw (q1q2) -- (q1) -- (q1q3) -- (q3) ;
    \draw (q1q2) -- (q2) -- (q2q3);
    \draw (q2q3) -- (q3);
\end{tikzpicture}

\subsection{Proof of Theorem \ref{item:11} and consequences}

Before proving Theorem \ref{item:11}, we have to show the following two lemmas. Lemma \ref{item:7} will allow us to prove Lemma \ref{boncardCK} - in which we prove $C(\K)$ has cardinality $r_1+r_2-1$, which is the first step in the proof of Theorem \ref{item:11}. For any $i \in \llbracket 1,r \rrbracket$, let $d_i$ denote the degree of $\K (i)/\mbb{Q}$.

\begin{lemme} \label{item:7}
    For any non-empty subset $\Omega \subset \llbracket 1 , r \rrbracket $, let
    \begin{gather*}
    f_{\mbb{C}} (\Omega) = \frac{1}{2} \prod\limits_{i \in \Omega} \left(  d_i-1 \right) + \frac{(-1)^{|\Omega|}}{2},  \quad f_{\mbb{R}} (\Omega) =  \prod\limits_{i \in \Omega} \left( d_i-1 \right)\\
        g_{\mbb{C}} (\Omega) = \frac{1}{2} \prod\limits_{i \in \Omega} d_i, \quad g_{\mbb{R}} (\Omega) = \prod\limits_{i \in \Omega} d_i .
    \end{gather*}
Let each of these functions map $\emptyset$ to $1$.
We have
    \begin{align}
        \sum\limits_{\substack{ \Omega \subset \llbracket 1 , r \rrbracket \\ \Omega \neq \emptyset }} f_{\mbb{C}}(\Omega) &= g_{\mbb{C}}(\llbracket 1 , r \rrbracket) -1 \label{item:16} \\
        \sum\limits_{\substack{\Omega \subset \llbracket 1 , r \rrbracket \\ \Omega \neq \emptyset }} f_{\mbb{R}}(\Omega) &= g_{\mbb{R}}(\llbracket 1 , r \rrbracket) -1. \label{item:17}
    \end{align}
\end{lemme}

\begin{proof}
    \underline{We start with Equation \eqref{item:17}.}

    We have to prove
    \[ \mbf{1} * f_{\mbb{R}}( \llbracket 1 , r \rrbracket) = g_{\mbb{R}}(\llbracket 1 , r \rrbracket)\]
    but instead, we will show that we have 
    \[ \forall \Omega \subset \llbracket 1 , r \rrbracket, \quad f_{\mbb{R}}( \Omega) = \mu *g_{\mbb{R}}(\Omega) \]
    and Equation \eqref{item:17} will follow from Theorem \ref{thmconvprod}.
    We have
    \begin{align*}
        \mu *g_{\mbb{R}}(\Omega) &= \sum\limits_{X \subset \Omega} (-1)^{|\Omega|-|X|} g_{\mbb{R}}(X) \\
        &= \sum\limits_{X \subset \Omega} (-1)^{|\Omega|-|X|} \prod\limits_{i \in X} d_i \\
        &= (-1)^{|\Omega|} + \sum\limits_{k=1}^{|\Omega|}  (-1)^{|\Omega|-k} \sum\limits_{\substack{i_1,\dots,i_k \in \,  \Omega \\ i_1 < \dots < i_k}} d_{i_1} \cdots d_{i_k}.
    \end{align*}

   Using Vieta's formulas, we see that this last expression matches the evaluation of the polynomial $(-1)^{|\Omega|} \prod\limits_{i \in \Omega} T-d_i$ at $T=1$, hence
    \[\mu *g_{\mbb{R}}(\Omega) = (-1)^{|\Omega|} \prod\limits_{i \in \Omega} 1-d_i = f_{\mbb{R}}(\Omega).\]
    In a similar way, \underline{we now consider Equation \eqref{item:16}.}
    We have
    \begin{align*}
        \mu *g_{\mbb{C}}(\Omega) &= \sum\limits_{X \subset \Omega} (-1)^{|\Omega|-|X|} g_{\mbb{C}}(X) \\
        &= (-1)^{\Omega} + \sum\limits_{\substack{X \subset \Omega \\ X \neq \emptyset}} (-1)^{|\Omega|-|X|} \frac{1}{2 }\prod\limits_{i \in X} d_i \\
        &= (-1)^{\Omega} + \underbrace{\frac{1}{2 }\sum\limits_{k=1}^{|\Omega|}  (-1)^{|\Omega|-k} \sum\limits_{\substack{i_1,\dots,i_k \in \, \Omega \\ i_1 < \dots < i_k}} d_{i_1} \cdots d_{i_k}}_{=\frac{(-1)^{|\Omega|}}{2} \left( \left(\prod\limits_{i \in \Omega} T-d_i \right) - T^{|\Omega|} \right)_{T=1}}
    \end{align*}
by Vieta's formulas again, hence
\[ \mu *g_{\mbb{C}}(\Omega) = (-1)^{\Omega} + \frac{(-1)^{|\Omega|}}{2} \left( \left(\prod\limits_{i \in \Omega} 1-d_i \right) - 1 \right) =f_{\mbb{C}}(\Omega). \]
\end{proof}

\begin{lemme} \label{boncardCK}
    The family $C(\mbb{K})$ has cardinality $r_1+r_2-1$.
\end{lemme}

\begin{proof}
    We have to show
    \[  \abs{C(\mbb{K})}  = \left\{
    \begin{array}{ll}
        \frac{1}{2} \left( \prod\limits_{i \in \llbracket 1 ,r \rrbracket} d_i \right) -1 & \mbox{if } \llbracket 1 ,r \rrbracket _{\mbb{C}} \neq \emptyset  \\
        \left( \prod\limits_{i \in \llbracket 1 ,r \rrbracket} d_i \right) -1 & \mbox{otherwise.}
    \end{array}
\right. \]
This can also be stated in the following way. For any non-empty subset $\Omega \subset \llbracket 1 , r \rrbracket$, let $f(\Omega)=\abs{C_{\Omega}(\mbb{K})}$ and let

\[g(\Omega) =\left\{
    \begin{array}{ll}
        \frac{1}{2} \prod\limits_{i \in \Omega} d_i & \mbox{if } \Omega_{\mbb{C}} \neq \emptyset  \\
        \prod\limits_{i \in \Omega} d_i & \mbox{otherwise.}
    \end{array}
\right. \]
Also, let $f$ and $g$ map $\emptyset$ to $1$. Then, we have to show
\[ \mbf{1}*f (\llbracket 1 ,r \rrbracket) = g(\llbracket 1 ,r \rrbracket).\]
Again, we rather show 
\begin{equation} \label{eqprodconv}
\forall \Omega \subset \llbracket 1 ,r \rrbracket, \quad f(\Omega) = \mu * g(\Omega)
\end{equation}
and Theorem \ref{thmconvprod} will conclude.
Let $\Omega \subset \llbracket 1,r \rrbracket$. If $ \abs{\Omega} \leqslant 1$, a straightforward computation shows Equation \eqref{eqprodconv} holds for $\Omega$. Suppose $\abs{\Omega} > 1$ and let $i_1 < \dots < i_s$ be such that $\Omega=\{i_1, \dots, i_s \}$.
We separate three cases.

\underline{Suppose we have $\Omega_{\mbb{C}} = \emptyset $.}
Then, Lemma \ref{item:7} gives
\[
\mu * g (\Omega) = \mu * g_{\R}(\Omega) = \prod\limits_{i \in \Omega} (d_i -1) 
\]
and it remains to observe this product equals $f(\Omega)$ as we supposed $\Omega_{\mbb{C}} = \emptyset $.

\underline{Now suppose $\Omega_{\mbb{R}} = \emptyset$.}
For any integer $k$, let $C_{\Omega}^k(\K)$ denote the elements of $C_{\Omega}$ that are of the form $u_1 \cdots u_k$.
If $\abs{\Omega}$ is odd, we have
\[ f(\Omega) =  \sum\limits_{k=1}^s | C_{\Omega}^k(\mbb{K}) | = \sum\limits_{k=1}^s (\frac{1}{2}d_{i_k} -1)(d_{i_{k-1}}-1) \cdots (d_{i_1} -1)  \]
and a straightforward induction on $l \in \llbracket 1,s \rrbracket $ shows
\begin{equation} \label{eqreccard}
\sum\limits_{k=1}^l (\frac{1}{2}d_{i_k} -1)(d_{i_{k-1}}-1) \cdots (d_{i_1} -1) = \frac{1}{2}(d_{i_l}-1) \cdots (d_{i_1}-1) -\frac{1}{2}.
\end{equation}
Taking $l =s$ and considering Lemma \ref{item:7}, we get
$ \mu * g (\Omega) = f_{\mbb{C}}(\Omega) = f(\Omega)$.
If $\abs{\Omega}$ is even, we have
\[ f(\Omega) = 1+ \sum\limits_{k=1}^s | C_{\Omega}^k(\mbb{K}) | = 1+ \sum\limits_{k=1}^s (\frac{1}{2}d_{i_k} -1)(d_{i_{k-1}}-1) \cdots (d_{i_1} -1)  \]
and we get the same conclusion.

\underline{Now, suppose that we have $\Omega_{\mbb{C}}  \neq  \emptyset$ and $\Omega_{\mbb{R}} \neq \emptyset$.} We have

\begin{align*}
    \mu * g (\Omega) =& \sum\limits_{X \subset \Omega} (-1)^{|\Omega| - |X|} g(X) \\
                =& \sum\limits_{\substack{X_1 \subset \Omega_{\mbb{R}} \\ X_2 \subset \Omega_{\mbb{C}} \\ X_2 \neq \emptyset}} (-1)^{|\Omega| - |X_1| - |X_2|} \times \frac{1}{2} \prod\limits_{i \in X_1 \cup X_2} d_i + \sum\limits_{X_1 \subset \Omega_{\mbb{R}}} (-1)^{|\Omega| - |X_1|} \prod\limits_{i \in X_1} d_i \\
                =& \sum\limits_{\substack{X_1 \subset \Omega_{\mbb{R}} \\ X_2 \subset \Omega_{\mbb{C}}}} (-1)^{|\Omega| - |X_1| - |X_2|} g_{\mbb{C}}(X_1 \cup X_2) \\
                & - \sum\limits_{X_1 \subset \Omega_{\mbb{R}}} (-1)^{|\Omega_{\mbb{C}}|  + |\Omega_{\mbb{R}}| - |X_1| } g_{\mbb{C}}(X_1) \\
                &+ \sum\limits_{X_1 \subset \Omega_{\mbb{R}}} (-1)^{|\Omega_{\mbb{C}}| + |\Omega_{\mbb{R}}| - |X_1|} g_{\mbb{R}}(X_1) \\
                =& f_{\mbb{C}}(\Omega) - (-1)^{|\Omega_{\mbb{C}}|} f_{\mbb{C}}(\Omega_{\mbb{R}}) + (-1)^{|\Omega_{\mbb{C}}|} f_{\mbb{R}}(\Omega_{\mbb{R}}) \\
                =& \frac{1}{2} \prod\limits_{i \in \Omega} (d_i-1) + \frac{(-1)^{|\Omega_{\mbb{C}}|}}{2} \prod\limits_{i \in \Omega_{\mbb{R}}} (d_i-1) .
\end{align*}   

Separate cases depending on whether $\abs{\Omega_{\mbb{C}}}$ is even or not and a similar induction argument to that of the proof of Equation \eqref{eqreccard} shows
\[ f(\Omega) = \frac{1}{2} \prod\limits_{i \in \Omega} (d_i-1) + \frac{(-1)^{|\Omega_{\mbb{C}}|}}{2} \prod\limits_{i \in \Omega_{\mbb{R}}} (d_i-1) .\]
\end{proof}

\begin{thm} \label{item:11}
The family $C(\mbb{K})$ is a $\Z[ 1/2 ]$-basis of $\was_2(\mbb{K})$.
Moreover, the $\Z[ 1/2 ]$-module $\was_2(\mbb{K})$ is a direct factor of $\was_2(\mbb{Q}(\zeta_n))$.
\end{thm}

\begin{proof}
Let us show the elements of $C(\mbb{K})$ generate a direct factor of $\was_2(\mbb{Q}(\zeta_n))$.
To this aim, we will investigate the decomposition in the basis $B$ (see Theorem \ref{theoremkucbase}) of every element of $C(\K)$. More precisely, for any non-empty subset $\Omega \subset \llbracket 1,r \rrbracket$, for any element of $C_{\Omega}(\K)$, we will investigate the part of their decomposition that lies in $B_{\Omega}$. If $\abs{\Omega} \geqslant 2$, we will associate to each $c \in C_{\Omega}(\K)$ a term $\phi(c) \in B_{\Omega}$ such that $\phi(c)$ appears with exponent $1$ or $2$ in the decomposition of $c$ in the basis $B$. Moreover, given distinct elements $c_1,c_2 \in C_{\Omega}(\K)$, we will show
\begin{enumerate}[(P1)] \hypertarget{P1}{}
    \item \label{P1} $\phi(c_1)$ is not involved in the decomposition of $c_2$ if both $c_1$ and $c_2$ are not problematic terms 
    \item \label{P2} $\phi(c_1)$ is not involved in the decomposition of $c_2$ if $c_1$ is a problematic term.
\end{enumerate}
If $\abs{\Omega}=1$, we will also associate to each $c \in C_{\Omega}(\K)$ a term $\phi(c) \in B_{\Omega}$ but the situation is a bit different and will be explained later. In particular, in this case, we will show $\phi(c)$ appears with exponent $-1$ in the decomposition of $c$ in $B$.

We will make use of \hyperlink{P1}{(P1)} and \hyperlink{P1}{(P2)} to order the elements of $B$ and also order the terms of $C(\K) \cup (B \setminus \phi(C(\K)))$ so that the matrix of this last family in the basis $B$ is triangular with diagonal coefficients lying in $\{  \pm 1, 2 \}$. This matrix is thus invertible in $\mbb{Z} [1/2]$ and so $C(\K)$ generates a direct factor of $\was_2(\Q(\zeta_n))$. Then, as Lemma \ref{boncardCK} gives the cardinality of $C(\K)$, we can apply Lemma \ref{item:18} to conclude. To ease the reading, we will handle elements of $C_{\llbracket 1 ,r \rrbracket}(\K)$ only but it is clear that the same kind of arguments works for any other $C_{\Omega}(\K)$ (those results concerning any $C_{\Omega}(\K)$ can also be obtained as a consequence of the case $\Omega={\llbracket 1 ,r \rrbracket}$ since we have $C_{\Omega}(\K)=C_{\Omega}(\K_{\Omega})$).

Recall $b_n$ is defined in Definition \ref{defcn}. Let $u \in Y_{\llbracket 1 ,r \rrbracket}(\K)$ and let $c = c_{n}(\K)^u$. We will show that we can let 
\[
\phi(c) = \left\{
    \begin{array}{ll}
        b_{n}^u & \mbox{ if } r \geqslant 2 \\
        \xi_{n,u} & \mbox{ otherwise}. 
    \end{array}
\right. 
\] In each of the following cases, we will then compute the exponent of $b_{n}^u$ in the decomposition of $c$ and we will investigate the decomposition of $c$.

\underline{Suppose $r = 1$.}
Modulo roots of unity, we have
\begin{align} 
    c &= \prod\limits_{w \in \Gal(\Q(\zeta_{n})^{+}/\K^{+})} \xi_{n,\sigma_1}^{uw} \nonumber \\
    &=\prod\limits_{w \in \mc{T}_{1}(\K^{+})} \frac{1-\zeta_n^{\sigma_1 uw}}{1-\zeta_n^{uw}} \nonumber \\
    &=\prod\limits_{w \in \mc{T}_{1}(\K^{+})} \frac{(1-\zeta_n)^{\sigma_1 uw}}{1-\zeta_n} \frac{1-\zeta_n}{1-\zeta_n^{uw}} \nonumber \\
    c &= \prod\limits_{w \in \mc{T}_{1}(\K^{+})} \xi_{n,\sigma_1 uw} \xi_{n,uw}^{-1}. \label{casr1}
\end{align}
Note that we have modulo roots of unity
\begin{gather*}
    \forall w \in \Gal(\Q(\zeta_n)/\Q), \quad \xi_{n,w} = \xi_{n,Jw} \\
    \xi_{n,1} = 1 \\
    \forall w \in \Gal(\Q(\zeta_n)/\Q) \setminus \{1,J \}, \quad w \in X_{\{1\}}(\Q(\zeta_n)) \mbox{ or } Jw \in X_{\{1\}}(\Q(\zeta_n))
\end{gather*}
and plugging these facts into Equation \eqref{casr1} gives the decomposition of $c$ in $B$. Observe there is a unique $(u',w) \in \mc{R}_{1,1}(\K^{+}) \times \mc{T}_{1}(\K^{+})$ such that $\sigma_1 u = u' w \mod J$ and define $\sigma_1 * u = u'$. Then, for any $c' \in C(\K)$, we see $\phi(c')$ appears in the decomposition of $c$ if and only if $c'= c$ or $c'=c_n(\K)^{\sigma_1 * u}$, with exponent $-1$ and $1$ respectively (we have $\sigma_1 uw \neq uw' \mod J $ for any $w,w' \in \mc{T}_1(\K^+)$, otherwise we would have $\sigma_1 \in \Gal(\Q(\zeta_n)/\K^{+})$, that is $\K^{+} = \Q$, so that $C(\K)$ is actually empty). Be aware that we may have $c_n(\K)^{\sigma_1 * u} \not\in C(\K)$ (this happens only when $\sigma_1 * u = J_1$) and this will explain why the top right coefficient of some matrix is not $1$ later. Define the following sequence
\begin{gather*}
    u^{(1)} = J_1 \\
    \forall i \in \llbracket 2, [ \K^{+} : \Q]  \rrbracket, \quad u^{(i)} = \sigma_1 * u^{(i-1)}
\end{gather*}
where $\sigma_1 * J_1$ is defined similarly as $\sigma * u$ and note that we have
\[
\left\{ u^{(i)} : i \in \llbracket 1, [ \K^{+} : \Q]  \rrbracket \right\} = \mc{R}_{1,1}(\K^{+})
\]
because $\sigma_1^{0}, \dots, \sigma_{1}^{[ \K^{+} : \Q] -1}$ is also a set of representatives of $\Gal(\K^{+}/\Q)$.
Define a strict total order $<_{\{1\},1}$ on $Y_{\{1\}}(\K)$ by setting
\begin{equation} \label{defordre}
    \forall i,j \in \llbracket 2, [ \K^{+} : \Q]  \rrbracket, \quad u^{(i)} <_{\{1\},1} u^{(j)} \Longleftrightarrow{} i < j.
\end{equation}

\underline{From now on, suppose $r \geqslant 2$.}

\underline{Suppose $\llbracket 1 ,r \rrbracket_{\C} = \emptyset$.}
This case has already been considered in \cite[Proposition 2, Remark 4]{werl} and we now write it down for convenience.

Modulo roots of unity of $\mbb{Q}(\zeta_n)$, we have
\begin{align} \label{decomporeel}
    c&= \op{N}_{\mbb{L}_n / \mbb{K}} (\eta_n^{u}) \nonumber \\
    &= \prod\limits_{w_1 \in \mc{T}_1(\K)} \cdots \prod\limits_{w_r \in \mc{T}_r(\K)} \prod\limits_{\varepsilon_1,\dots,\varepsilon_{r-1} \in \left\{ 0;1 \right\}} 1-\zeta_n^{J_1^{\varepsilon_1} u_1 w_1 \cdots J_{r-1}^{\varepsilon_{r-1}}u_{r-1}w_{r-1} u_r w_r }
\end{align}
and this is the decomposition of $c$ in the basis $B$. Indeed, for any $i < r$, we have $ J_i^{\varepsilon_i}u_{i}w_i \in \Gal(\Q(\zeta_{q_i}) / \Q) \setminus \{ J_i \}$ and $u_r w_r \in \mc{R}_r \setminus \{ 1 \}$ by construction of those $\mc{R}_{i,1}(\K),\mc{T}_i(\K)$.

As expected, we see $\phi(c) = 1-\zeta_n^u$ appears with exponent 1. Indeed, by construction of those $\mc{R}_{i,1}(\K),\mc{T}_i(\K)$, the products of the form $J_1^{\varepsilon_1} u_1 w_1 \cdots J_{r-1}^{\varepsilon_{r-1}}u_{r-1}w_{r-1} u_r w_r$ that appear in Equation \eqref{decomporeel} are pairwise distinct. Also, in this case, observe the decomposition of any $c' \in C_{\llbracket 1, r \rrbracket}(\K) \setminus \{ c \} $ is disjoint from the decomposition of $c$ as any $1-\zeta_n^w$ that appears in the decomposition of $c$ satisfies $w=u \mod \Gal(\Q(\zeta_n)/\K)$. In particular, we have \hyperlink{P1}{(P1)}.

\underline{Suppose $\llbracket 1 ,r \rrbracket_{\R} = \emptyset$.} Suppose $2 \nmid n$ (we will explain what to do if $2 \mid n$ later). We go through two cases depending on $u$ (and the parity of $r$).

Suppose $u \neq 1$. Then, we have
\begin{equation} \label{eqzetau}
    N_{\Q(\zeta_n)/\K} \left( 1- \zeta_n^u \right) = \prod\limits_{ s_1 \in \Gal(\Q(\zeta_{q_1})/ \K (1)) } \cdots \prod\limits_{s_r \in \Gal(\Q(\zeta_{q_r})/ \K (r))} 1-\zeta_n^{u s_1 \cdots s_r}
\end{equation}
and observe that we have 
\[
\forall \ s_1 \in \Gal(\Q(\zeta_{q_1})/ \K (1)), \dots, s_r \in \Gal(\Q(\zeta_{q_r})/ \K (r)), \quad u s_1 \cdots s_r \in X_{\llbracket 1,r \rrbracket},
\]
so that Equation \eqref{eqzetau} is the decomposition of $c$ in $B$.
Then, we see $\phi(c) = 1-\zeta_n^{u}$ appears with exponent $1$ in the decomposition of $c$ (because, again, those $u s_1 \cdots s_r$'s are pairwise distinct by construction of those $\mc{R}_{i,1}(\K),\mc{R}_{i,2}(\K)$). Note that any $1-\zeta_n^w$ that appears in this decomposition satisfies $w=u$ modulo $\Gal(\Q(\zeta_n) / \K)$ so that the decomposition of any $c' \in C_{\llbracket 1, r \rrbracket}(\K) \setminus \{ c, c_n(\K) \} $ is disjoint from that of $c$.

Suppose $u=1$ (this case has to be considered when $r$ is even only). We have the same Equation \eqref{eqzetau} as before and the same observations can be made for the same reasons (we have $1-\zeta_n \in B_{\llbracket 1,r \rrbracket}$ as $r$ is even). More precisely, the decomposition of any $c' \in C_{\llbracket 1,r \rrbracket}(\K) \setminus \{ c \}$ is disjoint from that of $c$ and $\phi(c)$ appears with exponent $1$ in the decomposition of $c$. In particular, we have \hyperlink{P1}{(P1)} and \hyperlink{P1}{(P2)}.

If $2 \mid n$, we have to do more manipulations to get the decomposition of $c$. First, recall we suppose $p_r=2$ in this case. Write $u=u_1 \cdots u_k $ with $k \in \llbracket 0,r \rrbracket$, $u_i \in \mc{R}_{i,1}(\K) \setminus \{ J_i \}$ for any $i \in \llbracket 1, k-1 \rrbracket$ and $u_k \in \mc{R}_{k,2}(\K) \setminus \{1 \}$. Then, let $u_i = 1$ for any $i > k$ so that $u=u_1 \cdots u_r$. We still have Equation \eqref{eqzetau} and now consider the notation introduced in this equation. If $u_r s_r \in \mc{R}_r$, then we still have $1-\zeta_n^{u s_1 \cdots s_r} \in B$ - in particular, this happens when $s_r=1$. If $u_rs_r \not\in \mc{R}_{r}$, we can show that $1-\zeta_n^{u s_1 \cdots s_r}$ decomposes in $B$ with
\begin{itemize}
    \item terms of $B$ whose level is lower than $n$,
    \item terms of the form $1-\zeta_n^{ s_{1}'  \cdots s_{r-1}' s_{r} '}$ with $s_r'=J_r u_r s_r \in \mc{R}_r$ and for any $i \in \llbracket 1 , r-1 \rrbracket$, $s_i' \in \Gal(\Q(\zeta_{q_i})/\Q) \setminus \{ J_i \}$.
\end{itemize}
Indeed, Equation \eqref{conj} gives
\[
 1-\zeta_n^{u s_1 \cdots s_r} = 1-\zeta_n^{u J_1 s_1 \cdots J_r s_r}.
\]
 Then, Lemma \ref{keyf1} shows that, modulo roots of unity, the element $1-\zeta_n^{u J_1s_1 \cdots J_rs_r}$ is a product of 
 \begin{itemize}
    \item elements of $B$ whose level is lower than $n$
    \item elements of the form
 \[
 1-\zeta_n^{ \pm s_{1}'  \cdots s_{r-1}' s_{r} '}
 \]
 such that $s_i' \in \Gal(\Q(\zeta_{q_i})/\Q) \setminus \{ J_i \}$ and $s_r'=J_r u_r s_r$. Note that we have $ 1-\zeta_n^{ s_{1}'  \cdots s_{r-1}' s_{r} '} \in B$ and $ 1-\zeta_n^{ s_{1}'  \cdots s_{r-1}' s_{r} '} \neq 1-\zeta_n^{u}$.
 \end{itemize}
 Finally, we conclude again that $\phi(c)=1-\zeta_n^{u}$ appears with exponent $1$ but, this time, the decompositions of the elements of $C_{\llbracket 1,r \rrbracket}(\K)$ may not be pairwise disjoint. We still conclude that $\phi(c)$ is not involved in the decomposition of any element of $C_{\llbracket 1,r \rrbracket}(\K) \setminus \{ c \}$ as we just showed the following. If $w=w_1 \cdots w_j$ for some $j \in \llbracket 0,r \rrbracket$ is involved in the decomposition of $c$, let $w_{j+1} =1, \dots, w_r = 1$; then we have one of the following two cases
\begin{enumerate}[i)]
    \item $w=u$ modulo $\Gal(\Q(\zeta_n)/\K)$
    \item $w_r \neq u_r$ and $w_r = u_r$ modulo $\Gal(\Q(\zeta_{q_r})/\K (r)^{+})$.
\end{enumerate}
Regardless of the parity of $n$, if some $w \neq u$ is such that $1-\zeta_n^w$ appears in the decomposition of $c$, then we have $w \not\in Y_{\llbracket 1,r \rrbracket}(\K)$. In particular, we have \hyperlink{P1}{(P1)}.

\underline{Suppose $\llbracket 1 ,r \rrbracket_{\C} \neq \emptyset$ and $\llbracket 1 ,r \rrbracket_{\R} \neq \emptyset$.}
Suppose $u$ is \underline{not} of the form $u_1 \cdots u_{t-1}$ with $u_i \in \mc{R}_{i,1}(\K) \setminus \{ J_i \}$ (that is $u$ is not a problematic term).
We still have the same Equation \eqref{eqzetau} and similar statements can be made. More precisely, we see $\phi(c)$ appears with exponent $1$ in the decomposition of $c$. If $2 \nmid n$ or $\K\cap \Q(\zeta_{2^{v_2(n)}})$ is real, then, the decomposition of any non problematic $c' \in C_{\llbracket 1,r \rrbracket}(\K) \setminus \{ c \}$ is disjoint from the decomposition of $c$ and if some $1-\zeta_n^w$ appears in the decomposition of $c$, then we have $w=u \mod \Gal(\Q(\zeta_n) / \K)$. If $2\mid n$, if $\K\cap \Q(\zeta_{2^{v_2(n)}})$ is imaginary and some $1-\zeta_n^w$ appears in the decomposition of $c$, then we have one of the following cases:
\begin{enumerate}[i)]
    \item $w=u \mod\Gal(\Q(\zeta_n)/\K)$
    \item $w_r \neq u_r$ and $w_r = u_r \mod \Gal(\Q(\zeta_{q_r})/\K (r)^{+})$.
\end{enumerate}
Regardless of the parity of $n$, observe that we have $w \not\in Y_{\llbracket 1,r \rrbracket}(\K)$ whenever $w \neq u$ if $1-\zeta_n^w$ appears in the decomposition of $c$.  In particular, we have \hyperlink{P1}{(P1)}.

Now, suppose $u$ is problematic, that is $u=u_1 \cdots u_{t-1}$ with $u_i \in \mc{R}_{i,1}(\K) \setminus \{ J_i \}$ for any $i \in \llbracket 1, t-1 \rrbracket$ (this case has to be considered when $ \abs{\llbracket 1,r \rrbracket_{\C}}$ is even only). Recall that we defined a length in Lemma \ref{keyf2}. We will show $c$ decomposes in $B$ with
\begin{itemize}
    \item elements of $B_{\llbracket 1,r \rrbracket}$ whose length is greater than $t-1$,
    \item elements $1-\zeta_n^w$ with $w=w_1 \cdots w_{t-1}$ satisfying $
w=u \mod \Gal(\Q(\zeta_n)/\K)
$
    \item elements of $B$ whose level is lower than $n$
\end{itemize}
and $\phi(c)=1-\zeta_n^u$ appears with exponent $2$. Notice that, in the second bullet point, we have $w \not\in Y_{\llbracket 1,r \rrbracket}(\K)$ whenever $w \neq u$. Again, we have Equation \eqref{eqzetau}. Assuming $2 \nmid n$ or $\K \cap \Q(\zeta_{2^{v_2(n)}})$ is real, if one of the $s_i$'s is non trivial for some $i \in \llbracket t, r \rrbracket$, then we have $1-\zeta_n^{u s_1 \cdots s_r} \in B$. If $2\mid n$ and $\K \cap \Q(\zeta_{2^{v_2(n)}})$ is imaginary, we may not have $1-\zeta_n^{u s_1 \cdots s_r} \in B$ given $s_i \neq 1$ for some $i \in \llbracket t,r \rrbracket$ but Lemma \ref{keyf2} shows such $1-\zeta_n^{u s_1 \cdots s_r}$ still decomposes with elements of $B_{\llbracket 1,r \rrbracket}$ whose length is greater than $t-1$ and elements of $B$ whose level is lower than $n$. Then, regardless of the parity of $n$, we now just have to consider the decomposition of the following product
\[
\prod\limits_{ s_1 \in \Gal(\Q(\zeta_{q_1})/ \K (1)) } \cdots \prod\limits_{s_{t-1} \in \Gal(\Q(\zeta_{q_{t-1}})/ \K (t-1))} 1-\zeta_n^{u s_1 \cdots s_{t-1}}.
\]
For now, let $s_1 \in \Gal(\Q(\zeta_{q_1}) / \K (1)), \dots, s_{t-1} \in \Gal(\Q(\zeta_{q_{t-1}})/ \K (t-1))$. If we have $u_{t-1}s_{t-1} \in \mc{R}_{t-1}$ (that is if $s_{t-1} \in \mc{T}_{t-1}(\K)$), then $1-\zeta_n^{u s_1 \cdots s_{t-1}} \in B$. Under this condition, observe $1-\zeta_n^u$ appears if and only if $s_1=1, \dots, s_{t-1}=1$ and it appears with exponent $1$. Else (that is if $s_{t-1} \in J_{t-1}\mc{T}_{t-1}(\K)$), observe that we have $u_i s_i \neq 1$ for any $i \in \llbracket 1, t \llbracket$ so that Lemma \ref{keyf2} shows that $1-\zeta_n^{u s_1 \cdots s_{t-1}}$ decomposes with
\begin{itemize}
    \item elements of $B_{\llbracket 1,r \rrbracket}$ whose length is greater than $t-1$,
    \item elements of $B$ whose level is lower than $n$,
    \item $1-\zeta_n^{u J_1s_1 \cdots J_{t-1}s_{t-1}}$ and its exponent is $(-1)^{r-t+1}=1$.
\end{itemize}
Then, under the condition $u_{t-1}s_{t-1} \not\in \mc{R}_{t-1}$, we see $1-\zeta_n^u$ appears if and only if $s_1 = J_1, \dots, s_{t-1}=J_{t-1}$ and it appears with exponent $1$. In particular, we have \hyperlink{P1}{(P2)} and $\phi(c)$ appears with exponent $2$ if $c$ is problematic.

We are then done with decomposing the elements of $C(\K)$.
We can now construct the matrix we talked about earlier.

To this aim, define a strict total order $<_{lex}$ - that is the lexicographic order defined in \cite[page 19]{lexord} - on the powerset of $\llbracket 1,r \rrbracket$ as follows:
\[
    \forall \ \Omega_1,\Omega_2 \subset \llbracket 1,r \rrbracket, \quad \Omega_1 <_{lex} \Omega_2 \Longleftrightarrow \renewcommand{\arraystretch}{1.5} \left\{
        \begin{array}{*1{>{\displaystyle}l}}
            \abs{\Omega_1} < \abs{\Omega_2}   \mbox{ or } \\
            \abs{\Omega_1} = \abs{\Omega_2} \mbox{ and} \\
            \quad \min \Omega_1 \setminus (\Omega_1 \cap \Omega_2) <\min \Omega_2 \setminus (\Omega_1 \cap \Omega_2).
        \end{array}
           \right.
\]
For any non-empty set $\Omega = \{ i_1, \dots, i_s \} \subset \llbracket 1,r \rrbracket$, for any $k \in \llbracket t_{\Omega}-1, s \rrbracket$, let $Y_{\Omega}^k(\K)$ be the set made of all the elements of $Y_{\Omega}(\K)$ of the form $u_1 \cdots u_k$. If $\abs{\Omega} \geqslant 2$, let $<_{\Omega,k}$ be any strict total order on $Y_{\Omega}^k(\K)$. If $\abs{\Omega}=1$, let $<_{\Omega,1}$ denote the strict total order on $Y_{\Omega}(\K)$ that is analogous to the one defined by Equation \eqref{defordre} for $\Omega=\{1 \}$.

Now, we are about to construct a strict total order on $C(\K)$. Given non-empty subsets $\Omega_1,\Omega_2 \subset \llbracket 1, r \rrbracket$ and $u_1 \cdots u_{k_1} \in Y_{\Omega_1}(\K), w_1 \cdots w_{k_2} \in Y_{\Omega_2}(\K)$, say that we have $c_{\Omega_1}(\K)^{ u_1 \cdots u_{k_1}} <  c_{\Omega_2}(\K)^{w_1 \cdots w_{k_2}}$ if one the following three conditions is satisfied
\begin{gather*}
          \abs{\Omega_1} >_{lex} \abs{\Omega_2}  \\
        \Omega_1 = \Omega_2 \mbox{ and } k_1 < k_2 \\
        \Omega_1 = \Omega_2 \mbox{ and } k_1 = k_2 \mbox{ and } u_1 \cdots u_{k_1} <_{\Omega,k_1} w_1 \cdots w_{k_2} .
\end{gather*}
Then, there are $c^{(1)}, \dots, c^{({\abs{C(\K)}})} \in C(\K)$ such that $c^{(1)} < \dots < c^{(\abs{C(\K)})}$. Also, let $b^{(\abs{C(\K)}+1)}, \dots, b^{(\abs{C})} \in B \setminus \phi(C(\K))$ be pairwise distinct. Consider the matrix of $(c^{(1)}, \dots, c^{(\abs{C(\K)})},  b^{(\abs{C(\K)}+1)}, \dots, b^{(\abs{C})})$ in $(\phi(c^{(1)}), \dots, \phi(c^{(\abs{C(\K)})}), b^{(\abs{C(\K)}+1)}, \dots, b^{(\abs{C})})$ (and observe this last family is $B$).
We have a triangular matrix that is just as expected and we now explain why. First, by the first condition that defines our order on $C(\K)$, the matrix we constructed is of the following form
\[
    \left( 
        \begin{array}{cc}
            \begin{array}{ccc}
                M_{\llbracket 1,r \rrbracket} & & \\
                & \ddots & \\
                  &  & M_{\{ 1 \}}
            \end{array} & \makebox(-2em,2em){\text{\huge0}} \\ 
             \makebox(-1.5em,1em){\text{\huge*}} & I
        \end{array}
        \right)
\]
where 
\begin{itemize}
    \item $I$ denotes the identity matrix with size $\op{Card}(B \setminus \phi(C(\K)))$ (and represents the terms of $B \setminus \phi(C(\K))$)
    \item and each matrix $M_{\Omega}$ represents partially the $B_{\Omega}$-part of the decomposition of the elements of $C_{\Omega}(\K)$ in the basis $B$.
\end{itemize}
Note that the zeros appear as any term from $C_{\Omega}(\K)$ decomposes in $B$ with terms whose level is lower than or equal to $\Omega$ (as explained after Definition \ref{deflevel}).

Then, let $\Omega$ be a non-empty subset of $\llbracket 1,r \rrbracket$. If $\Omega$ has cardinality $1$, then $M_{\Omega}$ is of the following form
\[
\left(
    \begin{array}{cccc}
        -1 & & &  \makebox(-1em,-1em){\text{\huge0}} \\
        1 &\ddots  \\
        & \ddots & \ddots \\
        \makebox(0em,1em){\text{\huge0}} & & 1 & -1
    \end{array}
\right)
\]
as a result of the remarks we have made on the decomposition of $c_{n}(\K)^u$ in the case $r=1$ and the fact that our order takes account of those remarks. One may wonder why the top right coefficient of $M_{\Omega}$ is not $1$ in this case and this is explained by the fact that we have $\xi_{n_{\Omega},1} \not\in \phi(C_{\Omega}(\K))$ (see our explanations in the case $r=1$ above).

 From now, suppose $\abs{\Omega} \geqslant 2$. If $\Omega_{\C}=\emptyset$ or $\Omega_{\R} = \emptyset$, then $M_{\Omega}$ is the identity matrix because of \hyperlink{P1}{(P1)}. Now, suppose we have $\Omega_{\R} \neq \emptyset$ and $\Omega_{\C} \neq \emptyset$. The matrix $M_{\Omega}$ is the identity matrix if $\abs{\Omega_{\C}}$ is odd as shown by \hyperlink{P1}{(P1)}. If $\abs{\Omega_{\C}}$ is even, the matrix $M_{\Omega}$ is of the following form
\[
    \left( 
        \begin{array}{c|c}
            2I & 0 \\
            \hline * & I
        \end{array}
        \right)
\]
where the scaling matrix on the top left side corresponds to the problematic terms and the identity matrix on the bottom right side corresponds to the other terms of $C_{\Omega}(\K)$. Indeed, the second condition that defines our order on $C(\K)$ explains why the problematic terms appear first. Our fact \hyperlink{P1}{(P2)} explains why those zeros appear and our fact \hyperlink{P1}{(P1)} explains the identity matrix on the bottom right side.
Then, it suffices to use Lemma \ref{item:18} to conclude.
\end{proof}





This next corollary already results from \cite[Theorem 3.1]{KuceraRevue} but we still show how it is a consequence our $\Z[1/2]$-basis $C(\K)$.

\begin{coro} \rightskip 0pt plus 1fill
    Suppose $\K$ is totally deployed. The quotient group $\was(\mbb{K}) / \Sin(\mbb{K})$ is a $2$-group.
\end{coro}

\begin{proof} 
    Indeed, any $c_{\Omega}(\K)^{u} \in C(\K)$ is already an element of $\Sin(\mbb{K}_{\Omega})$ if $ \abs{\Omega_{\mbb{C}}} \geqslant 1 $ and $\abs{\Omega} \geqslant 2$ . Any other element of $C(\K)$ has order $1$ or $2$ in the quotient group $\was(\mbb{K}_{\Omega}) / \Sin(\mbb{K}_{\Omega})$ (see \cite[Equation 11, Corollary 3]{werl}).
    Hence, the quotient group $(\was(\mbb{K}) / \Sin(\mbb{K}) ) \otimes \Z [ 1 /2 ]$ is trivial. 
\end{proof}

    If $\K (1),\dots,\K (r)$ are real, Werl Milàn stated and proved in a special case (see \cite[Remark 4]{werl}) this quotient group is an elementary $2$-group with rank $ [ \mbb{K}: \mbb{Q} ] -1 $ and the family $C(\K)$ is a $\Z$-basis of $\was(\K)$.

    If $\K (1),\dots,\K (r)$ are imaginary, observe we have $\was(\mbb{K}) = \Sin(\mbb{K})$ as shows the proof of Theorem \ref{item:11} (again, in this case, no $2$'s appear on the diagonal of the triangular matrix of this proof). Then $C(\K)$ is a $\Z$-basis of $\Sin(\mbb{K})$. The same $\Z$-basis of $\Sin(\mbb{K})$ has been given in \cite[Corollary 4.3]{KuceraBase} and the author also proved $\Sin(\K) = \was(\K)$ (it results from \cite[Theorem 5.1]{KuceraBase}).

    The other cases are covered by the following Corollaries \ref{alpha} and \ref{diffwassin}.

\begin{coro} \label{alpha}
    Suppose $\K$ is totally deployed. Let $M( \K)$ denote the abelian group generated by $C(\K)$ and $\mbf{Z}(\K)$. Assume one of the $\K(i)$'s is imaginary. We have 
    \[
    [\was(\K) : M( \K)] \leqslant 2^{\alpha}
    \]
    with 
    \[
    \alpha= (2^{r-t}-1) \left( [\K (1) \cdots \K (t-1) : \Q] -1 \right)  .
    \]
\end{coro}

\begin{proof}
    Quotient by the roots of unity and keep the same notation for $M(\K)$ and $\was$.
    
    Through the proof of Theorem \ref{item:11}, we see that we have a subgroup $T$ of $\was(\Q(\zeta_n))$ such that $T \cap M(\K) = \{ 1 \}$ and 
    \[
    [ \was(\Q(\zeta_n)) : T \oplus M(\K)] = 2^{\alpha}
    \]
    with
    \[
        \alpha= \sum\limits_{\substack{\Omega \subset \llbracket 1,r \rrbracket \\ \abs{\Omega_{\C}} \in 2 \N^{*} \\ \Omega_{\R} \neq \emptyset }} \prod\limits_{i \in \Omega_{\R}} (d_i-1).
    \]
    A straightforward computation shows this definition of $\alpha$ matches the value given in the statement of Corollary \ref{alpha}. 
Then, we have
    \[
    [\was(\K) \oplus T : M(\K) \oplus T ] \leqslant 2^{\alpha}.
    \]
    Now, observe the natural map $\was(\K) / M(\K) \to \was(\K) \oplus T / M(\K) \oplus T$ is injective. Indeed, let $x \in \was(\K)$ be such that $x=yz$ with $y \in M(\K)$ and $z \in T$. As $M(\K)$ has finite index in $\was(\K)$, there is $k \in \N$ such that $x^k \in M(\K)$. Then, we have $z^k \in T \cap M(\K)$ so that $z^k = 1 = z$ since $T$ is torsion-free and $x=y \in M(\K)$.
\end{proof}

\begin{coro} \label{diffwassin}
    We have $[\was(\K) : \Sin(\K) ] \leqslant 2^{2^{r-t}([\K(1) \cdots \K(t-1) : \Q]-1)}$, assuming $\K$ is totally deployed and one of the $\K(i)$'s is imaginary.
\end{coro}

\begin{prv}
     We have $ C(\K) \setminus \Sin(\K) \subset C_{\llbracket 1, t -1 \rrbracket}(\K)$ and the elements of $C_{\llbracket 1, t -1 \rrbracket}(\K)$ all have order at most $2$ in $\was(\K)/ \Sin(\K)$ (see \cite[Equation 11, Corollary 3]{werl}). We have $\abs{C_{\llbracket 1, t -1 \rrbracket}(\K)} = \abs{C(\K_{\llbracket 1,t-1 \rrbracket})} = [\K(1) \cdots \K(t-1) : \Q]-1$ by Lemma \ref{boncardCK} so that it suffices to call Corollary \ref{alpha} to conclude.
\end{prv}

\begin{coro} \label{item:14}
    Let $Q$ denote the Hasse's unit index of $\K$.
    Assuming $\K$ is totally deployed and one of the $\K(i)$'s is imaginary, we have
    \[ [\mbf{E}(\mbb{K}): \was(\mbb{K}) ] = h^{+}(\mbb{K}) Q 2^x \]
    for some $x \in \mbb{Z}$ satisfying
    \[
    -{2^{r-t}([\K(1) \cdots \K(t-1) : \Q]-1)} - \nu \leqslant x \leqslant - \mu
    \]
    with $\nu$ being the number of integers $i$ such that $\K (i) / \Q$ has even degree and $\mu=r-t+1$ being the number of integers $i$ such that $\K (i) $ is imaginary.
\end{coro}

\begin{proof} 
    This results from Corollary \ref{diffwassin} and the formula Sinnott has given for the index of $\Sin(\mbb{K})$ in $\mbf{E}(\mbb{K})$ (see \cite[Proposition 4.1, Theorem 4.1, Theorem 5.4]{Sinnott1980}).
\end{proof}

If $\K(1), \dots, \K(r)$ are real, we have $[\mbf{E}(\mbb{K}): \was(\mbb{K}) ] = [\mbf{E}(\mbb{K}): M(\mbb{K}) ] = h(\K)$ as explained in \cite[Remark 4]{werl}, where $M(\K)$ denotes the group generated by $C(\K)$ and $\mbf{Z}(\K)$.

This next corollary can also be obtained using class field theory but it still arises naturally from our $\Z[1/2]$-basis $C(\K)$.

\begin{coro} \label{item:19}
    Suppose $\K = \K (1) \cdots \K (r)$ is totally deployed.
    Let $(A_1,\dots,A_k)$ be a partition of $\llbracket 1 ,r \rrbracket$. We have a canonical injective map
    \[
    \prod\limits_{j=1}^k \mbf{E}(\mbb{K}_{A_j}) / \was (\mbb{K}_{A_j}) \otimes \Z [1/2] \xhookrightarrow{} \mbf{E} (\mbb{K}) / \was(\mbb{K}) \otimes \Z [1/2].
    \] 
In particular, if we let $h_p^{+}(\mbb{K})$ denote the $p$-part of the class number of $\mbb{K}^{+}$, we have for any odd prime $p$
    \[ \prod\limits_{j=1}^k h_p^{+}(\mbb{K}_{A_j}) \mid h_p^{+}(\mbb{K}).  \]
\end{coro}

\begin{proof}
    Let $x=x_1 \cdots x_{k} \in \was_2(\mbb{K})$ with $x_j \in \mbf{E}(\mbb{K}_{A_j}) \otimes \Z [1/2]$ for any $j \in \llbracket 1,k \rrbracket$. We have to show $x_j \in \was_2(\mbb{K}_{A_j})$ for any $j$. There is an integer $N$ such that $x_j^N \in \was_2(\mbb{K}_{A_j})$ for any $j$ as $\mbf{E}$ and $\was$ have the same rank. Modulo roots of unity of $\mbb{K}$, we can decompose $x$ in $C(\K)$ and $x_j^{N}$ in $C(\K_{A_j})$ for any $j$:
    \[
        x = \prod\limits_{c \in C(\mbb{K})} c^{a_{x,c}}, \quad x_j^N =  \prod\limits_{c \in C(\mbb{K}_{A_j} )} c^{a_{x_{j},c}}.
    \]
    The construction of $C(\K)$ gives $C(\K_{A_1}) \sqcup \dots \sqcup C(\K_{A_k}) \subset C(\K)$.
    Then, we also get the decomposition of $x^N$ in $C(\K)$ by injecting the decomposition of $x_1^N,\dots, x_k^N$ in $x=x_1 \cdots x_{k}.$, so that
    \[ \forall j \in \llbracket 1 , k \rrbracket, \ \forall c \in C(\mbb{K}_{A_j}), \quad Na_{x,c} = a_{x_{j},c} \]
    then we have
    \[ x_j^N = \left( \prod\limits_{c \in C(\mbb{K}_{A_j} )} c^{a_{x,c}} \right)^N \]
    hence $x_j \in \was_2(\mbb{K}_{A_j})$. The result on class numbers is then given by Corollary \ref{item:14} and \cite[Remark 4]{werl}.
\end{proof}

\section{Along the cyclotomic tower} \label{sectioninfini}

Through this section, let $p$ be an odd prime number. We will give a $\Lambda$-basis of
$\was_{\infty}(\K) = \lim\limits_{\xleftarrow{}} \left(\was(\K_k) / \mbf{Z}(\K_k) \otimes \Z_p \right)$
where
\begin{itemize}
    \item $\K$ is a totally deployed abelian number field
    \item $(\K_k)_{k \in \N}$ denotes the cyclotomic $\Z_p$-tower of $\K$
    \item $\Lambda$ denotes the Iwasawa algebra, that is $\Lambda=\lim\limits_{\xleftarrow{}} \Z_p[\Gal(\K_k/\K)]$ (the limit is taken with respect to the restriction maps)
    \item the limit defining $\was_{\infty}$ is taken with respect to the norm maps.
\end{itemize}
We will have two cases to handle for technical reasons. Our plan, in each of these cases, is to give a family $B^{\infty}(\K)$ of elements of $\was_{\infty}(\K)$ that generates $\was_{\infty}(\K)$ as a $\Lambda$-module and that is made of $r_1 + r_2$ elements (where $r_1$ is the number of real embeddings of $\K$ and $r_2$ is half of the number of complex embeddings of $\K$); this is enough to show that $B^{\infty}(\K)$ is a basis as it has been proven in \cite{bel} that $\was_{\infty}(\K)$ is a free $\Lambda$-module of rank $r_1+r_2$ (it results from \cite[Proposition 1.3, Corollary 1.6]{bel}). We will construct $B^{\infty}(\K)$ as follows. For each $k \geqslant 1$, we will first give a $\Z_p$-basis of $\was(\K_k)$ that we will denote by $C'(\K_k)$ (and that is a slight modification of the family $C(\K_k)$ that we gave in Section \ref{sectionnotCK}). As a consequence of Lemma \ref{pregeninf}, in order to get our $\Lambda$-basis $B^{\infty}(\K)$, we will know that we may simply make use of those terms of $C'(\K_k)$ that involve $1-\zeta_{\Omega}$ for some $\Omega$ such that $1 \in \Omega$; we can naturally construct elements of $\was_{\infty}(\K)$ from those special terms of the $C'(\K_k)$'s: that is what we will do with the map $\mc{T}$ (see Definition \ref{transfoT}). Then, we simply consider the action of $\Gal(\K_{\infty}/ \K)$ on the image of $\mc{T}$ and we give $B^{\infty}(\K)$ as a set of representatives of this action.

\subsection{Notation and preliminaries}

Let $\mbb{K} = \K{(1)} \cdots \K{(r)}$ be a totally deployed abelian number field of conductor $n$. If $\mbb{K}$ is not ramified at $p$, write $\K = \K{(2)} \cdots \K{(r)}, n=q_2 \cdots q_r$ instead and let $p_1=p$.

Let $\mbb{K} = \K_0 \subset \K_1 \subset \dots $ denote the cyclotomic $\Z_p$-tower of $\mbb{K}$ and let $\K_{\infty}$ denote the cyclotomic $\Z_p$-extension of $\K$ (that is the union of the $\K_k$'s). Define the Iwasawa's algebra associated to $\Gal(\K_{\infty}/\K)$ by $\Lambda=\lim\limits_{\xleftarrow{}} \Z_p[\Gal(\K_k/\K)]$ with respect to the restriction maps. Let $\overline{\was}$ denote $\was / \mbf{Z} \otimes \Z_p$ and let $\was_{\infty}(\K) = \lim\limits_{\xleftarrow{}} \overline{\was}(\K_k)$ with respect to the norm maps. It appears clearly that $\was_{\infty}(\K)$ is a $\Lambda$-module. Let $\widetilde{\was}_k$ denote the universal norms of $\was(\K_k)$, that is the image of the canonical morphism $\was_{\infty}(\K) \to \overline{\was}(\K_k)$.

Let $\K^{tame}$ denote the maximal subfield of $\K$ that is tamely ramified at $p$. As $\K$ is totally deployed, the field $\K^{tame}$ is also totally deployed. Note that we have a canonical isomorphism of $\Gal(\K_{\infty}/\Q)$-modules $\was_{\infty}(\K) \simeq \was_{\infty}(\K^{tame})$ so that we may suppose $\K$ is tamely ramified at $p$ with no loss of generality (then, note that $\K$ may be unramified at $p$).

For any prime number $l$, let $\Q(\zeta_{l^{\infty}})$ denote the compositum of the $\Q(\zeta_{l^m})$'s for $m$ running over $\N$. For any supernatural integer $m=\prod_i l_i^{m_i}$, let $\Q(\zeta_m)$ denote the compositum of the $\Q(\zeta_{l_i^{m_i}})$'s.

Let $n_1 = \prod\limits_{j=1}^r p_j^{e_j'} = \prod\limits_{j=1}^r q_{j,1}$ be the conductor of $\K_1$. Let $q_{j,k}$ be the $p_j$-part of the conductor of $\K_k$. For any non-empty set $\Omega \subset \llbracket 1,r \rrbracket$, let $\mbb{Q}(\zeta_{\Omega}^{\infty})$ denote the compositum of the $\mbb{Q}(\zeta_{p_i^{\infty}})$'s where $i$ runs over $\Omega$. For any $k \in \N^{*} \cup \{ \infty \}$, let
\[ 
\K_{k,\Omega} = \K_k \cap \mbb{Q}(\zeta_{\Omega}^{\infty}).
\]
Let $i_0$ be such that $p_{i_0} = p$. Now, let $\sigma_{i_0} \in \Gal(\Q(\zeta_{p^{\infty}})/\Q)$ be such that, for any $k$, the restriction of $\sigma_{i_0}$ generates $\Gal(\Q(\zeta_{p^{k}})/\Q)$ (so that there is no conflict between the definitions of the $C(\K_k)$'s).

If $j \neq i_0$ (resp. $j=i_0$), see the elements of $\Gal(\K_j/\mbb{Q})$ (resp. $\Gal(\K_{\infty,\{i_0\}}/\mbb{Q})$) as elements of $\Gal(\K_{\infty}/\mbb{Q})$ by letting them act trivially on the compositum of the $\K_{\infty,\{i\}}$'s for $i$ running over $\llbracket 1,r \rrbracket \setminus \{ j \}$.

Similarly, if $j \neq i_0$ (resp. $j=i_0$), let $J_j$ be the complex conjugation of $\Gal(\mbb{Q}(\zeta_{q_j})/\mbb{Q})$ (resp. $\Gal(\mbb{Q}(\zeta_{p^{\infty}})/\mbb{Q})$) seen as an element of $\Gal(\mbb{Q}(\zeta_{np^{\infty}})/\mbb{Q})$. Let $J= J_1 \cdots J_r$ denote the complex conjugation of $\Gal(\mbb{Q}(\zeta_{np^{\infty}})/\mbb{Q})$.

For any non-empty subset $\Omega \subset \llbracket 1, r \rrbracket$, let
\begin{gather*}
    c_{\Omega}'(\K) = \renewcommand{\arraystretch}{1.5}\left\{
    \begin{array}{*2{>{\displaystyle}l}}
       \op{N}_{\mbb{Q}(\zeta_{\Omega}) / \mbb{K}_{\Omega}} (\xi_{q_i,\sigma_i}) & \mbox{if } \abs{\Omega} = 1 \\
        \op{N}_{\mbb{Q}(\zeta_{\Omega}) / \mbb{K}_{\Omega}}  (1-\zeta_{\Omega}) & \mbox{if } \abs{\Omega}  > 1 .
    \end{array}
\right. \\
C_{\Omega}'(\K) = \left\{  c_{\Omega}'(\K)^u : u \in Y_{\Omega}(\K) \right\}.
\end{gather*}
Let $C'(\K) = \cup_{\Omega} C_{\Omega}'(\K)$ and observe $C'(\K)$ is a $\Z_p$-basis of $\overline{\was}(\K)$.
Indeed, Theorem \ref{item:11} implies that $C(\K)$ is a $\Z_p$-basis of $\overline{\was}(\K)$ and we have just squared those elements of $C(\K)$ that are associated to any $\Omega$ such that $\K_{\Omega}$ is real.

\begin{lemme} \label{pregeninf}
    Let $k \geqslant 1$. Then $\widetilde{\was}_k$ is generated as a $\Gal(\K_{k}^{(i_0)}/\mbb{Q})$-module by
\[ \bigcup\limits_{\substack{\Omega \subset \llbracket 1, r \rrbracket \\ i_0 \in \Omega}} C_{\Omega}'(\K_k) .
\]

\end{lemme}

\begin{prv}
    The idea is to consider $C'(\K_{(k+j)})$ as $j$ tends to infinity and apply the norm map from $\K_{(k+j)}$ to $\K_{(k)}$. Let $M$ denote the $\Gal(\K_k^{(i_0)}/\mbb{Q})$-module generated by
    \[ \bigcup\limits_{\substack{\Omega \subset \llbracket 1, r \rrbracket \\ i_0 \in \Omega}} C_{\Omega}'(\K_{k}) .\]
    Let $x \in \widetilde{\was}_k$. For any $j \in \mbb{N}$, there is $y \in \overline{\was}(\K_{(k+j)})$ such that $ x = \op{N}_{\K_{k+j} / \K_{k}} (y)$.
    We can decompose $y$ in the basis $C'(\K_{k+j})$, that is there are some $\alpha_{\Omega}$'s being a sum of terms of the form $a u$ with $a \in \mbb{Z}_p$ and $u \in Y_{\Omega}(\K_{k + j})$ such that
    \[
     y = \prod\limits_{\substack{\Omega \subset \llbracket 1,r \rrbracket \\ \Omega \neq \emptyset }} c_{\Omega}'(\K_{k+j})^{\alpha_{\Omega}}.
     \]
    Hence we have
    \begin{align*}
        x &= \left( \prod\limits_{\substack{\Omega \subset \llbracket 1,r \rrbracket \\ i_0 \in \Omega }} \op{N}_{\K_{k+j} / \K_{k}} (c_{\Omega}'(\K_{k+j})^{\alpha_{\Omega}}) \right) \left( \prod\limits_{\substack{\Omega \subset \llbracket 1,r \rrbracket \setminus \{ i_0 \} \\ \Omega \neq \emptyset }} \op{N}_{\K_{k+j} / \K_k} (c_{\Omega}'(\K_{k+j})^{\alpha_{\Omega}}) \right) \\
        &= \left( \prod\limits_{\substack{\Omega \subset \llbracket 1,r \rrbracket \\ i_0 \in \Omega }}  c_{\Omega}'(\K_{k})^{\alpha_{\Omega}} \right) \left( \prod\limits_{\substack{\Omega \subset \llbracket 1,r \rrbracket \setminus \{ i_0 \} \\ \Omega \neq \emptyset }}  c_{\Omega}'(\K_{k+j})^{p^{j}\alpha_{\Omega}} \right).
    \end{align*} 

If $i_0 \not\in \Omega$, we have $c_{\Omega}'(\K_{k+j}) = c_{\Omega}'(\mbb{K}_1)$ (which does not depend on $j$).
Then, we have written $x=m_j {m_j'}^{p^{j}}$ with $m_j \in M$ and $m_j' \in \overline{\was}(\mbb{K}_1)$. As both $M$ and $\overline{\was}(\mbb{K}_1)$ can be seen as finitely generated free $\mbb{Z}_p$-module, we see that $M$ is compact (for the $p$-adic topology) and we have ${m_j'}^{p^{j}} \xrightarrow[j \to \infty]{} 1$. Then, taking the limit in $x=m_j {m_j'}^{p^{j}}$ along a converging subsequence of $(m_j)_{j \geqslant 0}$ shows that we have $x \in M$.
\end{prv}




\begin{definition} \label{transfoT}
    Let $\mc{T}$ be the following transformation: given any $\Omega \subset \llbracket 1,r \rrbracket$ such that $ i_0 \in \Omega$ and any $u \in \Gal(\K_{{\infty},\Omega} / \Q)$, define $\mc{T}(c_{\Omega}'(\K_{1})^{u}) \in \was_{\infty}(\mbb{K})$ as the sequence whose projection on $\widetilde{\was_k}$ for any $k \geqslant 1$ is $c_{\Omega}'(\K_{k})^{u}$.
\end{definition}

Note that we should write $\mc{T}(c_{\Omega}'(\K_{1}),u)$ instead of $\mc{T}(c_{\Omega}'(\K_{1})^{u})$ but we will keep this last notation.

\subsection{Unramified or real}

Through this subsection, suppose either $\mbb{K}$ is not ramified at $p$ or $\mbb{K}$ is tamely ramified at $p$ and $\K{(i_0)}$ is real. Then, for convenience, suppose $i_0=1$ (that is to say $p_1=p$). Recall $d_1$ denotes the degree of $\K{(1)} / \Q$ and, if $\mbb{K}$ is not ramified at $p$, let $d_{1} = 1$ and $\K{(1)}=\mbb{Q}$.

Let ${\Omega} = \{ 1 \} \subset \llbracket 1 , r \rrbracket $ and let 
\[
     X_{\Omega}^{\infty}(\mbb{K}) = \left\{ \sigma_1^{a} : a \in \llbracket 0, d_1 \llbracket \right\} .
\]
The restriction to $\K{(1)}$ induces a bijection between $X_{\{ 1 \}}^{\infty}(\mbb{K})$ and $\Gal(\K{(1)} / \Q)$. In fact, we could replace $X_{\{ 1 \}}^{\infty}(\mbb{K})$ with any lift of $\Gal(\K{(1)} / \Q)$ to $\Gal(\K_{\infty,\{1\}}/\Q)$ and we chose this one for convenience.

Recall, up to reordering, there is an integer $t$ such that $\K_1, \dots, \K_{t-1}$ are real and $\K_t, \dots, \K_r$ are imaginary. To understand more easily the following notation, note that, for any $i$, if $\K_i$ is real then we have $J_i=1$ when we see $J_i $ as an element of $\Gal(\K_i/\Q)$.

    For any ${\Omega} = \{ i_1,\dots,i_s \} \subset \llbracket 1 , r \rrbracket $ with $1 \in \Omega, s \geqslant 2$ and $i_1 < \dots < i_s$ and $\K{(i_s)}$ is imaginary (that is $\mbb{K}_{\Omega}$ decomposes with at least one imaginary field), recall $t_{\Omega}$ is the integer such that $\K{(i_1)},\dots,\K{(i_{t_{\Omega}-1})}$ are real and $ \K{(i_{t_{\Omega}})}, \dots, \K{(i_s)}$ are imaginary.
    Let $X_{\Omega}^{\infty}(\K)$ be the set of products of the form $u_1 \cdots u_k$ with $k \in \llbracket t_{\Omega},s \rrbracket$ and
    \[ u_1 \in X_{\{1\}}^{\infty}(\K), \quad \forall j \in \llbracket 2, k \llbracket, \ u_j \in \Gal(\K{(i_j)}/\Q) \setminus \{J_{i_j}\}, \quad u_k \in \mc{R}_{i_k,2}(\K) \setminus \{ 1 \}. \]
     If $\abs{\Omega_{\C}}$ is even, then add to $X_{\Omega}^{\infty}(\K)$ all the products of the form $u_1 \cdots u_{t_{\Omega}-1}$ with
     \[
     u_1 \in X_{\{1\}}^{\infty}(\K), \quad \forall j \in \llbracket 2, t_{\Omega }\llbracket, \ u_j \in \Gal(\K{(i_j)}/\Q) \setminus \{1\} .
     \]
If $\K{(i_{s})}$ is real (that is $\mbb{K}_{\Omega}$ decomposes with real fields only), let
\[
     X_{\Omega}^{\infty}(\mbb{K}) = \left\{ u_1 \cdots u_s : u_1 \in X_{\{1\}}^{\infty}(\K), \forall j > 1, \quad u_j \in \Gal(\K{(i_j)}/\Q) \setminus \{ 1 \} \right\}
\]
    For any non-empty subset $\Omega \subset \llbracket 1,r \rrbracket$ such that $1 \in \Omega$, let
\[
    B_{\Omega}^{\infty}(\K) = \mc{T} \left\{ c_{\Omega}'(\K_{1})^{u} : u \in X_{\Omega}^{\infty} \right\}.
    \]

Let $B^{\infty}(\mbb{K})=\cup_{\Omega} B_{\Omega}^{\infty}(\mbb{K})$ where $\Omega$ runs over the set of all subsets of $\llbracket 1 , r \rrbracket$ that contain $1$.

\begin{lemme} \label{cardBinf1}
    The family $B^{\infty}(\mbb{K})$ has cardinality $r_1 + r_2$.
\end{lemme}

\begin{prv}
     Let $d = r_1 + r_2$. Let $D$ be $1+ \op{rank}_{\Z}\was(\K{(2)} \cdots \K{(r)})$ and observe $D$ does not depend on $p$. \underline{First, suppose $\mbb{K}$ is not ramified at $p$}. For any non-empty $\Omega \subset \llbracket 1,r \rrbracket$, let $f(\Omega)$ denote the cardinality of $C_{\Omega}'(\K_{1})$. Now, from Lemma \ref{boncardCK}, applied to both $\K_{1}$ and $\K$, we get
    \begin{align*}
         \sum\limits_{\substack{\Omega \subset \llbracket 1 ,r \rrbracket \\ 1 \in \Omega}} f(\Omega) &= \sum\limits_{\substack{\Omega \subset \llbracket 1 ,r \rrbracket \\ \Omega \neq \emptyset}} f(\Omega) - \sum\limits_{\substack{\Omega \subset \llbracket 2 ,r \rrbracket \\ \Omega \neq \emptyset}} f(\Omega) \\
        &= p D -1 -(D-1) \\
        &= (p-1)D.
    \end{align*}
     Note that $f(\Omega)$ is a polynomial in variable $p$ with coefficients in $\Z$. As the last equality remains true for infinitely many $p$ (more precisely, for any $p$ that is distinct from the other $p_i's$), the last equality can be interpreted in terms of polynomials. Now, \underline{stop assuming $\mbb{K}$ is not ramified at $p$}. Evaluate at $d_1+1$ this last equality between polynomials to obtain the expected result as $d=d_1 D$ since $\K{(1)}$ is real (the swap between $p$ and $d_1+1$ is explained by the fact that, when $\K$ is unramified at $p$, $\abs{\Gal(\K_{1,\{1\}}/\Q) \setminus \{ 1 \}} = p-1$ controls the cardinality of $C_{\Omega}'(\K_{1})$ and that of $B^{\infty}(\K)$ is controlled by $\abs{\Gal(\K{(1)} / \Q)}=d_1$ in general).
\end{prv}

\begin{thm} \label{thmBinf1}
    The family $B^{\infty}(\mbb{K})$ is a $\Lambda$-basis of $\was_{\infty}(\mbb{K})$.
\end{thm}

\begin{prv}
    Lemma \ref{cardBinf1} implies that it only remains to show $B^{\infty}(\mbb{K})$ generates $\was_{\infty}(\mbb{K})$.
    To this aim, for any $k \geqslant 1$, we will prove that the projection of $B^{\infty}(\mbb{K})$ onto $\widetilde{\was_k}$ generates $\widetilde{\was_k}$ as a $\Z_p[\Gal(\K_k/\K)]$-module.
    
    Let $k \geqslant 1$ and let $x \in \widetilde{\was_k}$. Observe $\Gal(\K_{k} / \mbb{K})$ is generated by $\sigma_{1}^{d_1}$ and this gives $\Gal(\K_{k} / \Q) = X_{\{ 1 \}}^{\infty}(\mbb{K})\Gal(\K_{k} / \mbb{K})$. Hence, Lemma \ref{pregeninf} states $x$ lies in the $\Gal(\K_{k} / \mbb{K})$-module generated by the projection of $B^{\infty}(\K)$ onto $\widetilde{\was}_k$.
\end{prv}

\subsection{Imaginary and ramified}

Through this subsection, suppose $\mbb{K}$ is tamely ramified at $p$ and $\K{(i_0)}$ is imaginary. Then, suppose $i_0=t$ to ease the definition of our future basis $B^{\infty}(\K)$.
  
    Let ${\Omega} = \{ t \} \subset \llbracket 1 , r \rrbracket $ and let
    \[
        X_{\Omega}^{\infty}(\mbb{K}) = \mc{R}_{t,2}(\K) \subset \Gal(\K_{\infty} / \Q).
    \]   
One must understand the embedding $\mc{R}_{t,2}(\K) \subset \Gal(\K_{\infty} / \Q)$ through the fact that we naturally have $\Gal(\K_{\infty}/\Q) \simeq \Gal(\K / \Q) \times \Gal(\K_{\infty}/\K)$ as $\K$ is assumed to be tamely ramified at $p$.

    Let $\Omega = \{ i_1, \dots , i_s \} \subset \llbracket 1,r \rrbracket$ with $t \in \Omega, s \geqslant 2$ and $i_1 < \dots < i_s$. Let $t_{\Omega}$ be such that $i_{t_{\Omega}} = t$. Let $X_{\Omega}^{\infty}(\mbb{K})$ be the set of products of the form $u_1 \cdots u_k$ with $k \in \rrbracket t_{\Omega},s \rrbracket$, satisfying $u_k \in \mc{R}_{i_k,2}(\K) \setminus \{ 1 \}$, $u_{t_{\Omega}} \in \Gal(\K{(t)}/ \Q) \subset \Gal(\K_{\infty} / \Q)$ and 
    \[
        \forall j \in \llbracket 1,k-1 \rrbracket \setminus \{ t_{\Omega} \}, \ u_j \in \Gal(\K{(i_j)}/\Q) \setminus \{ J_{i_j} \}.
    \]
    Add to $X_{\Omega}^{\infty}(\mbb{K})$ all the products $u_1 \cdots u_{t_{\Omega}}$ with
    \[
     u_{t_{\Omega}} \in \mc{R}_{t,2}(\K) \subset \Gal(\K_{\infty} / \Q), \quad \forall j \in \llbracket 1, t_{\Omega} \llbracket, \ u_j \in \Gal(\K{(i_j)}/\Q) \setminus \{ 1 \}.
    \]
    For any $\Omega \subset \llbracket 1,r \rrbracket$ such that $t \in \Omega$, let
    \[
        B_{\Omega}^{\infty}(\K) = \mc{T} \left\{ c_{\Omega}'(\K_{1})^{u} : u \in X_{\Omega}^{\infty}(\mbb{K}) \right\}.
    \]
Let $B^{\infty}(\mbb{K})=\cup_{\Omega} B_{\Omega}^{\infty}(\mbb{K})$ where $\Omega$ runs over the set of all subsets of $\llbracket 1 , r \rrbracket$ that contains $t$.

\begin{lemme} \label{cardBinf2}
    The family $B^{\infty}(\mbb{K})$ has cardinality $r_1+r_2$.
\end{lemme}

\begin{prv}
    Let $d=r_1+r_2=r_2$. Let $f(\Omega)$ denote the cardinality of $B_{\Omega}^{\infty}(\K)$ and, for any $i \in \llbracket 1 ,r \rrbracket$, recall $d_i = [ \K{(i)} : \Q ]$. We must prove
    \[ 
    \sum\limits_{\substack{\Omega \subset \llbracket 1,r \rrbracket \\ t \in \Omega} } f(\Omega) = \frac{d_1 \cdots d_r}{2}.
    \]
   Let \[g(\Omega) = \frac{1}{2}\prod\limits_{i \in \Omega} d_i .\]

   Let $E=\llbracket 1,r \rrbracket \setminus \{ t \}$. We want to show that we have
   \[
   \sum\limits_{\Omega \subset E} f(\Omega \cup \{ t \}) = g(E \cup \{ t \}).
   \]
   To this aim, we will prove that we have for any $\Omega \subset E$
   \[
   f(\Omega \cup \{ t \} ) = \sum\limits_{X \subset \Omega} \mu(\Omega \setminus X) g(X \cup \{ t\}).
   \]
   (see Theorem \ref{thmconvprod}). Indeed, we have
   \begin{align*}
       \sum\limits_{X \subset \Omega} \mu(\Omega \setminus X) g(X \cup \{ t\}) &= \sum\limits_{k=0}^{\abs{\Omega}} (-1)^{\abs{\Omega}-k} \sum\limits_{i_1 < \dots < i_k \in \Omega} \frac{d_{i_1} \cdots d_{i_k} d_t}{2} \\
       &= \frac{(-1)^{\abs{\Omega}} d_t}{2} P(1) \\
       &= \frac{d_t}{2} \prod\limits_{i \in \Omega} \left( d_i-1 \right)
   \end{align*}

   where $P(T) = \prod\limits_{i \in \Omega} \left( T-d_i \right)$ (we gathered the $X$'s according to their cardinality before using Vieta's formula).

   To conclude, it suffices to show that we have 
   \[
   f(\Omega \cup \{t \}) = \frac{d_t}{2} \prod\limits_{i \in \Omega} \left( d_i-1 \right).
   \]
   Let $i_1 < \dots < i_s$ be such that $\Omega \cup \{t \} = \{ i_1, \dots, i_s \}$ and let $t_{\Omega}'$ be such that $i_{t_{\Omega}'} = t$.
   For any $k \in \llbracket t_{\Omega}', s\rrbracket$, let $f(\Omega,k)$ denote the number of products of the form $u_1 \cdots u_k$ lying in $X_{\Omega \cup \{ t \}}^{\infty}(\mbb{K})$. A straightforward induction on $l \in \llbracket t_{\Omega}', s \rrbracket$ shows that we have
\[
\sum\limits_{k=t_{\Omega}'}^l f(\Omega,k) = \frac{d_t}{2} \prod\limits_{j \in \llbracket 1,l \rrbracket \setminus \{ t_{\Omega}' \}} \left( d_{i_j}-1 \right)
\]
and taking $l=s$ concludes.
\end{prv}

\begin{thm} \label{thmBinf2}
    The family $B^{\infty}(\mbb{K})$ is a $\Lambda$-basis of $\was_{\infty}(\mbb{K})$.
\end{thm}

\begin{prv}
    Lemma \ref{cardBinf2} implies that it only remains to show $B^{\infty}(\mbb{K})$ generates $\was_{\infty}(\mbb{K})$. To this aim, for any $k \geqslant 1$, we will prove  that the projection of $B^{\infty}(\mbb{K})$ onto $\widetilde{\was_k}$ generates $\widetilde{\was_k}$ as a $\Z_p[\Gal(\K_k/\K)]$-module.

    Let $k \geqslant 1$ and let $\gamma$ be a topological generator of $\Gal(\K_{\infty,\{t\}}/\K(t)) \simeq \Z_p$. Lemma \ref{pregeninf} states it suffices to make sure the $\Gal(\K_{k}/ \K)$-module generated by the projection of $B^{\infty}(\mbb{K})$ onto $\widetilde{\was_k}$ reaches all the elements of 
    \[ \bigcup\limits_{\substack{\Omega \subset \llbracket 1, r \rrbracket \\ t \in \Omega}} C_{\Omega}'(\K_{k}) .
    \]
    To show that, we just need to observe $\Gal(\K_{k,\{t\}}/ \K)$ is generated by (the restriction of) $\gamma$ and has odd order $p^k$ so that $\mc{R}_{t,2}(\K_{k}) = \langle \gamma \rangle \mc{R}_{t,2}(\K)$ (again, we naturally see $\mc{R}_{t,2}(\K) \subset \mc{R}_{t,2}(\K_{k})$ because $\K$ is supposed to be tamely ramified at $p$).
\end{prv}

\textbf{Acknowledgments} The author acknowledges financial support from ANR project PadLEfAn.

\bibliographystyle{plainurl}
\bibliography{refs}
\let\thefootnote\relax\footnote{Rafik SOUANEF \\ Université Marie et Louis Pasteur, CNRS, LmB (UMR 6623)\\16 route de Gray, 25000, Besançon, France\\Email: rafik.souanef@ens-rennes.fr\\Url: https://perso.eleves.ens-rennes.fr/people/rafik.souanef/}

\end{document}